\documentclass[11pt,reqno]{amsart}
\usepackage[T1]{fontenc}
\usepackage{lmodern}
\usepackage[english]{babel}
\usepackage{amsmath}
\usepackage[OT2,OT1]{fontenc}
\usepackage{amssymb}
\usepackage{mathrsfs}
\usepackage{bookmark}
\usepackage[numbers,sort&compress]{natbib}
\usepackage{slashed}
\usepackage{cancel}
\usepackage{hyperref}
\usepackage{amsthm}
\usepackage[a4paper]{geometry}
\begin{document}

\title[The Cohomology of the Witt and the Virasoro algebra]{The Vanishing of the Low-Dimensional Cohomology of the Witt and the Virasoro algebra}
\author[Jill Ecker]{Jill Ecker}
\address{University of Luxembourg, Faculty of Science, Technology and Communication,
Campus Belval, Maison du Nombre, 6, avenue de la Fonte, L-4364 Esch-sur-Alzette, Luxembourg}
\curraddr{}
\email{jill.ecker@uni.lu}
\author[Martin Schlichenmaier]{Martin Schlichenmaier}
\address{University of Luxembourg, Faculty of Science, Technology and Communication,
Campus Belval, Maison du Nombre, 6, avenue de la Fonte, L-4364 Esch-sur-Alzette, Luxembourg}
\curraddr{}
\email{martin.schlichenmaier@uni.lu}
\thanks{Partial  support by the
Internal Research Project  GEOMQ15,  University of Luxembourg,
and
by the OPEN programme  of the Fonds National de la Recherche
(FNR), Luxembourg,  project QUANTMOD O13/570706
is gratefully acknowledged.}
\subjclass[2000]{Primary: 17B56; Secondary: 17B68, 17B65, 17B66, 14D15, 
81R10, 81T40}

\keywords{Witt algebra; Virasoro algebra; Lie algebra cohomology; Deformations of algebras; conformal field theory}
\date{v1: 19.07.2017; v2: 31.01.2018} 
\begin{abstract}
A proof of the vanishing of the third cohomology group of the Witt algebra with values in the adjoint module is given. 
Moreover, we provide a sketch of the proof of the one-dimensionality of the third cohomology group of the Virasoro algebra with values in the adjoint module.    
The proofs given in the present article are completely algebraic and independent of any underlying topology. They are a generalization of the ones provided by Schlichenmaier, who proved the vanishing of the second cohomology group of the Witt and the Virasoro algebra by using purely algebraic methods. In the case of the third cohomology group though, extra difficulties arise and the involved proofs are distinctly more complicated. The first cohomology group can easily be computed; we will give an explicit proof of its vanishing in the appendix, in order to illustrate our techniques.  
\end{abstract}

\maketitle

\section{Introduction}
The \textit{Witt algebra} $\mathcal{W}$ is an infinite-dimensional, $\mathbb{Z}$-graded Lie algebra first introduced by Cartan in 1909 \cite{MR1509105}. The Witt algebra and its universal central extension, the \textit{Virasoro algebra}, are two of the most important infinite-dimensional Lie algebras, used in mathematics as well as in theoretical physics, see e.g. the book by Kac, Raina and Rozhkovskaya \cite{MR3185361}. Therefore, knowledge of their cohomology groups is of outermost importance to a better understanding of central extensions, outer morphisms, deformations, obstructions, and so on.\\
In the present article, we consider algebraic cohomology mainly with values in the adjoint module. The vanishing of the second cohomology group of the Witt and the Virasoro algebra with values in the adjoint module has been proved by Schlichenmaier in \cite{MR3200354,MR3363999} by using elementary algebraic methods; see also Fialowski \cite{MR2985241}. In \cite{MR1054321}, Fialowski announced the result for the Witt algebra without proof.
In \cite{MR1417181} by van den Hijligenberg and Kotchetkov, the vanishing of the second cohomology group with values in the adjoint module of the superalgebras $k(1),k^+(1)$ and of their central extensions was proved using algebraic methods.
The aim of this article is to prove the vanishing of the third cohomology group with values in the adjoint module of the Witt algebra by merely algebraic means. We shall provide a generalization of the proof given in \cite{MR3200354,MR3363999} to the third cohomology group of the Witt algebra. As our goal is to consider arbitrary cohomology, we are dealing completely with algebraic cocycles, meaning we do not need to put any continuity constraints on them.\\
The results obtained for algebraic cohomology might be compared to those obtained for continuous cohomology.
In fact, the Witt algebra is related to the group of diffeomorphisms of the unit circle $S^1$. It forms a dense subalgebra of the complexified algebra $Vect(S^1)$ of all smooth vector fields on $S^1$. Based on results of Goncharova \cite{MR0402765}, Reshetnikov \cite{MR0292097} and Tsujishita \cite{MR0458517}, the continuous cohomology of the Lie algebra $Vect(S^1)$ with values in the adjoint module has been computed by Fialowski and Schlichenmaier in \cite{MR2030563}, yielding:
\begin{equation*}
\mathrm{H}^*(Vect(S^1),Vect(S^1))=\{0\}\,.
\end{equation*}
As a first naive guess, one might be inclined to argue, by using density arguments, that this result implies the vanishing of the cohomology of the Witt algebra. However, this reasoning is not sound, as the argument using density is only valid when considering the sub-complex formed by the continuous cohomology of the Lie algebra. Instead, a direct proof for the vanishing of the algebraic cohomology is needed and furnished here for the low-dimensional cohomology.\\
Although the proof of the vanishing of the third cohomology group with values in the adjoint module of the Witt algebra is rather involved, it consists of very simple, elementary algebraic manipulations. Only basic knowledge of the cohomology of Lie algebras is needed to understand the proofs in this article. The basic notions necessary to understand the proofs will be introduced in the following sections. \\
The article is organized as follows: in Section \ref{Section2}, we introduce the Witt and the Virasoro algebra.\\ 
In Section \ref{Section3}, we give some basic notions of the Chevalley-Eilenberg cohomology of Lie algebras, and internally graded Lie algebras in particular. This includes the interpretation of the first four cohomology groups with values in the adjoint and the trivial module of a Lie algebra. Moreover, we present an overview of the results known for these cohomology groups in the case of the Witt algebra, including the results derived in this article.\\
Section \ref{Section5} consists of the main result of this article: It contains the proof of the vanishing of the third cohomology group with values in the adjoint module of the Witt algebra, i.e. $\mathrm{H}^3(\mathcal{W},\mathcal{W})=\{0\}$. The provided proofs are rather detailed in order to augment their readability.  \\
The final section, Section \ref{Section6}, contains the results for the Virasoro algebra. More precisely, we start by giving a summary of the results known for the low-dimensional cohomology of the Virasoro algebra, including the results obtained in this article. Subsequently, we show that the third cohomology of the Virasoro algebra with values in the trivial and the adjoint module is one-dimensional. Simultaneously, one can prove that the third cohomology with values in the trivial module is one-dimensional also in the case of the Witt algebra. In order to avoid overloading the present article, we only provide a sketch of the proofs involving the third cohomology of the Virasoro algebra. Details will be given in forthcoming work \cite{Ecker}.\\
The appendix contains the proofs related to the first cohomology of the Witt and the Virasoro algebra. The first result proven is  the vanishing of the first cohomology group with values in the adjoint module of the Witt algebra, i.e. $\mathrm{H}^1(\mathcal{W},\mathcal{W})=\{0\}$. This result means that all derivations of the Witt algebra are inner, which is a well-known fact, see Zhu and Meng \cite{MR1775845}. Still, we prove this result again explicitly in our algebraic setting, as the proof of $\mathrm{H}^1(\mathcal{W},\mathcal{W})=\{0\}$ is a nice warm-up example of our techniques and procedures. 
The proof is realized via a direct computation similar to the one performed in Section \ref{Section5}. The second result proven is the vanishing of the first cohomology group with values in the adjoint module of the Virasoro algebra, i.e. $\mathrm{H}^1(\mathcal{V},\mathcal{V})=\{0\}$. This proof uses techniques based on long exact cohomology sequences.\\
In the present article we only consider the first four cohomology groups of the Witt and the Virasoro algebra. These cohomology groups have a nice interpretation in terms of important Lie algebra objects such as invariants, outer derivations, central extensions, deformations, obstructions and crossed modules, see e.g. Gerstenhaber \cite{MR0161898,MR0171807,MR0207793,MR0240167,MR0160807}. 
Contrary to the first four cohomology groups, higher cohomology groups $\mathrm{H}^k(\mathcal{W},\mathcal{W})$ $k> 3$ do not come with such an easy interpretation in terms of known objects.  
Nonetheless, it would be interesting to have a similar elementary analysis also for higher cohomology groups.\\
Unfortunately, purely algebraic proofs become increasingly complicated when considering higher cohomology groups. A direct comparison of the proofs provided for the first, second and third cohomology groups already suggests that a proof of the vanishing of e.g. $\mathrm{H}^4(\mathcal{W},\mathcal{W})$ would be distinctly more intricate. In fact, although there are parallels between the three proofs, there is no immediate generalization when going from one cohomology group to the next cohomology group in the cohomology complex. Instead, each time new difficulties arise which cannot be anticipated when considering the precedent cohomology group. Still, we expect $\mathrm{H}^k(\mathcal{W},\mathcal{W})=\{0\}$ for all $k\geq 0$. 
\section{The Witt and the Virasoro algebra}\label{Section2}
The Witt algebra is generated as a vector space over a field $\mathbb{K}$ with characteristic zero by elements $\{e_n\ |\ n\in\mathbb{Z}\}$ which satisfy the following Lie algebra structure equation:
\begin{equation}
\left[e_n,e_m\right]=(m-n)e_{n+m}\,, \qquad n,m\in\mathbb{Z} \,. \label{StructureEquation}
\end{equation}
By defining the degree of an element $e_n$ by $deg(e_n):=n$, we obtain that the Witt algebra is a $\mathbb{Z}$-graded Lie algebra. Hence, it can be decomposed into infinitely many one-dimensional homogeneous subspaces, i.e. $\mathcal{W}=\bigoplus_{n\in\mathbb{Z}}\mathcal{W}_n$, where each subspace $\mathcal{W}_n$ is generated as a vector space by one basis element $e_n$. The Witt algebra is an \textit{internally} graded Lie algebra, i.e. the grading is given by one of its own elements, namely $e_0\in\mathcal{W}$, as can be seen from the Witt structure equation (\ref{StructureEquation}) yielding $\left[e_0,e_n\right]=n\ e_n=deg(e_n)\ e_n$. Thus, $e_n$ is an eigenvector of the adjoint action $ad_{e_0}$ with eigenvalue $n$. The associated eigenspace corresponds to $\mathcal{W}_n$. In other words, the grading of the Witt algebra corresponds to the eigenspace decomposition of the $ad_{e_0}$-action on $\mathcal{W}$.\\
As the second cohomology group of the Witt algebra with values in the trivial module is one-dimensional, i.e. $dim(\mathrm{H}^2(\mathcal{W},\mathbb{K}))=1$, there is, up to equivalence and rescaling, exactly one non-trivial central extension, which is the \textit{Virasoro} algebra $\mathcal{V}$. As a vector space, the Virasoro algebra is given as a direct sum $\mathcal{V}=\mathbb{K}\oplus\mathcal{W}$ generated by the basis elements $\hat{e}_n:=(0,e_n)$, $n\in\mathbb{Z}$ and its one-dimensional center $t:=(1,0)$. These generators satisfy the following Lie structure equation:
\begin{equation*}
[\hat{e}_n,\hat{e}_m]=(m-n)\hat{e}_{n+m}-\frac{1}{12}(n^3-n)\delta_n^{-m}t,\qquad [\hat{e}_n,t]=[t,t]=0\,,
\end{equation*}
for all $n,m\in\mathbb{Z}$\  \footnote{The Kronecker delta $\delta_n^{m}$ is defined as being equal to $1$ if $n=m$, and zero else.}.  
By defining $deg(\hat{e}_n):=deg(e_n)=n$ and $deg(t)=0$, also the Virasoro algebra $\mathcal{V}$ becomes an (internally) $\mathbb{Z}$-graded Lie algebra.  

\section{The cohomology of Lie algebras}\label{Section3}

\subsection{The Chevalley-Eilenberg cohomology}
In this section, for the convenience of the reader, we will recall the Chevalley-Eilenberg cohomology \cite{MR0024908}, which is the cohomology of Lie algebras.\\
Let $\mathcal{L}$ be a Lie algebra and $M$ an $\mathcal{L}$-module.
Moreover, let $C^q(\mathcal{L},M)$ be the vector space of $q$-multilinear alternating maps with values in $M$,
\begin{equation*}
C^q(\mathcal{L},M):=\text{Hom}_\mathbb{K}(\wedge^q\mathcal{L},M)\,.
\end{equation*}
The elements of $C^q(\mathcal{L},M)$ are called $q$-\textit{cochains}. By convention, one sets $C^0(\mathcal{L},M):=M$. The coboundary operators $\delta_q$ are defined by:
\begin{align}
\begin{array}[h]{rl}
\forall q\in\mathbb{N},\qquad &\delta_q:C^q(\mathcal{L},M)\rightarrow C^{q+1}(\mathcal{L},M): \psi \mapsto \delta_q \psi \,,\\
&\\
 (\delta_q\psi)(x_1,\dots x_{q+1}):&=\sum_{1\leq i<j\leq q+1}(-1)^{i+j+1}\ \psi (\left[x_i,x_j\right],x_1,\dots , \hat{x}_i,\dots , \hat{x}_j,\dots ,x_{q+1})\\
&\\
& +\sum_{i=1}^{q+1}(-1)^i\ x_i\cdot \psi (x_1,\dots ,\hat{x}_i,\dots ,x_{q+1})\,,
\end{array}\label{Delta}
\end{align}
with $x_1,\dots , x_{q+1}\in\mathcal{L}$, $\hat{x}_i$ means that the entry $x_i$ is omitted and the dot $\cdot$ stands for the module structure. In the case of the adjoint module, the module $M$ corresponds to the Lie algebra $\mathcal{L}$ itself and the module structure corresponds to the Lie algebra structure, i.e. $\cdot=\left[\cdot ,\cdot\right]$. In case of the trivial module $\mathbb{K}$, this product is zero. 
Since the coboundary operators fulfill $\delta_{q+1}\circ\delta_q=0\ \forall\ q\in\mathbb{N}$, we obtain a cochain complex $(C^*(\mathcal{L},M),\delta)$ which is called the \textit{Chevalley-Eilenberg complex}. The associated cohomology is the \textit{Chevalley-Eilenberg cohomology} \cite{MR0024908} given by:
\begin{equation*} 
\mathrm{H}^q(\mathcal{L},M):=Z^q(\mathcal{L},M)/B^q(\mathcal{L},M)\,,
\end{equation*}
where $Z^q(\mathcal{L},M):=\text{ ker }\delta_q$ is the vector space of $q$-\textit{cocycles} and 
$B^q(\mathcal{L},M):=\text{ im }\delta_{q-1}$ is the vector space of $q$-\textit{coboundaries}. 

\subsection{Degree of homogeneous cochains}
In the case of graded Lie algebras such as the Witt and the Virasoro algebra, the notion of degree for cochains can be introduced as a helpful tool. This will be important in the following.  
Since the present article focuses mainly on cohomology with values in the adjoint module, we will concentrate on the case $M=\mathcal{L}$.\\
Let $\mathcal{L}$ be a $\mathbb{Z}$-graded Lie algebra $\mathcal{L}=\bigoplus_{n\in\mathbb{Z}}\mathcal{L}_n$. A $q$-cochain $\psi$ is called \textit{homogeneous of degree d} if there exists a $d\in\mathbb{Z}$ such that for all $q$-tuple $x_1,\dots ,x_q$ of homogeneous elements $x_{i}\in\mathcal{L}_{deg(x_{i})}$, we have:
\begin{equation*}
\psi(x_{1},\dots ,x_{q})\in\mathcal{L}_n \text{ with }n=\sum_{i=1}^q deg(x_{i})+d\,.
\end{equation*}
In fact, the cohomology can be decomposed as follows:
\begin{equation*}
\mathrm{H}^q(\mathcal{L},\mathcal{L})=\bigoplus_{d\in\mathbb{Z}}\mathrm{H}^q_{(d)}(\mathcal{L},\mathcal{L})\,,
\end{equation*} 
where $\mathrm{H}^q_{(d)}(\mathcal{L},\mathcal{L})$ is the subspace consisting of classes of $q$-cocycles of degree $d$ modulo $q$-coboundaries of degree $d$.
Thus, the inspection of the cohomology group $\mathrm{H}^q(\mathcal{L},\mathcal{L})$ can be performed by analyzing each of its components $\mathrm{H}^q_{(d)}(\mathcal{L},\mathcal{L})$ separately. \\
An import result states that in the case of internally $\mathbb{Z}$-graded Lie algebras, nonzero cohomology groups can  exist only for degree zero, see Fuks \cite{MR874337}:
\begin{subequations}\label{Fuks}
\begin{equation}
\mathrm{H}^q_{(d)}(\mathcal{L},\mathcal{L})=\{0\}\ \text{ for }\ d\neq 0\ ,
\end{equation}
\begin{equation}
\mathrm{H}^q(\mathcal{L},\mathcal{L})=\mathrm{H}^q_{(0)}(\mathcal{L},\mathcal{L})\ .
\end{equation}
\end{subequations}
Nevertheless, we prove this statement again explicitly in the case of the Witt algebra, since the corresponding proof presented in this article is short and very simple. The aim is to render this article as self-contained as possible. \\
The discussion above focused on the adjoint module, but it holds true for any internally graded module, in particular for the trivial module $\mathbb{K}=\bigoplus_{n\in\mathbb{Z}}\mathbb{K}_n$ with trivial grading $\mathbb{K}_0=\mathbb{K}$ and $\mathbb{K}_n=\{0\}$ for $n\neq 0$.
\subsection{Interpretation of the first four cohomology groups} 
The cohomology of Lie algebras has important interpretations in terms of the Lie structure of the Lie algebra and its deformations, see Gerstenhaber \cite{MR0161898,MR0171807,MR0207793,MR0240167}. As we will need this later, we recall various interpretations here, with focus on the trivial module $\mathbb{K}$ and the adjoint module $\mathcal{L}$. 
\subsubsection{Zeroth cohomology}\label{Zeroth cohomology} 
The zeroth cohomology $\mathrm{H}^0(\mathcal{L},M)$ corresponds to the space of $\mathcal{L}$-invariants of the module $M$, denoted by $^{\mathcal{L}}M$. In particular, we have:
\begin{equation*}
\mathrm{H}^0(\mathcal{L},\mathbb{K})=\mathbb{K}\qquad\text{  and }\qquad\mathrm{H}^0(\mathcal{L},\mathcal{L})=C(\mathcal{L})\,,
\end{equation*} 
 where $C(\mathcal{L})$ is the center of $\mathcal{L}$.
\subsubsection{First cohomology}\label{First cohomology}  
The first cohomology groups with values in the trivial module and the adjoint module are given by:
\begin{align*}
&\mathrm{H}^1(\mathcal{L},\mathbb{K})=\left(\frac{\mathcal{L}}{[\mathcal{L},\mathcal{L}]}\right)^*\\
\text{ and }& \mathrm{H}^1(\mathcal{L},\mathcal{L})=\frac{\text{Der}(\mathcal{L}) }{ad_\mathcal{L}}=\text{Out}(\mathcal{L})\,,
\end{align*}
where $*$ stands for the dual space, $\text{Der}(\mathcal{L})$ for the derivations for $\mathcal{L}$, $ad_\mathcal{L}$ for the inner derivations and $\text{Out}(\mathcal{L})$ for the outer derivations.
\subsubsection{Second cohomology}
The second cohomology group $\mathrm{H}^2(\mathcal{L},\mathbb{K})$ with values in the trivial module classifies central extensions of $\mathcal{L}$ up to equivalence.\\
The second cohomology group $\mathrm{H}^2(\mathcal{L},\mathcal{L})$ with values in the adjoint module classifies infinitesimal deformations of $\mathcal{L}$ modulo equivalent deformations. In particular, $\mathrm{H}^2(\mathcal{L},\mathcal{L})=\{0\}$ implies that the Lie algebra $\mathcal{L}$ is infinitesimally and formally rigid, see Fialowski and Fuchs \cite{MR1670210}, Fialowski \cite{MR981622,MR1054321}, Gerstenhaber \cite{MR0171807,MR0207793,MR0240167}, and Nijenhuis and Richardson \cite{MR0195995}. 
 \subsubsection{Third cohomology}
A sensible question to ask is whether a given infinitesimal deformation can be lifted to a formal deformation. In general, this is not always possible because obstructions appear, which are elements of $\mathrm{H}^3(\mathcal{L},\mathcal{L})$. For more details, see e.g. \cite{MR2381782}. \\
However, the third cohomology comes also with a more constructive point of view.
Indeed, the third cohomology has an interpretation in terms of crossed modules. More precisely, the third cohomology group $\mathrm{H}^3(\mathcal{L},M)$ of the Lie algebra $\mathcal{L}$ with values in the module $M$ classifies equivalence classes of crossed modules associated to $\mathcal{L}$ and $M$, see Wagemann \cite{MR2229486} and Gerstenhaber \cite{MR0207793,MR0160807}.  
\subsection{Results for the cohomology of the Witt algebra}\label{ResultsWitt}
In this section, we briefly present known results of the cohomology of the Witt algebra with values in the trivial module and the adjoint module, and we also present a brief summary of the new results obtained in the present article concerning this topic.\\
Obviously, from Section \ref{Zeroth cohomology}, we see that the following is true:
\begin{equation*}
\mathrm{H}^0(\mathcal{W},\mathbb{K})=\mathbb{K}\qquad\text{ and }\qquad \mathrm{H}^0(\mathcal{W},\mathcal{W})=\{0\}\,.
\end{equation*}  
Furthermore, using the facts that the Witt algebra $\mathcal{W}$ is a perfect Lie algebra and that all derivations of the Witt algebra $\mathcal{W}$ are inner, we obtain from Section \ref{First cohomology} the following results:
\begin{equation*}
\mathrm{H}^1(\mathcal{W},\mathbb{K})=\{0\}\qquad\text{ and }\qquad \mathrm{H}^1(\mathcal{W},\mathcal{W})=\{0\}\,.
\end{equation*} 
It is a well-known fact that all derivations of the Witt algebra are inner derivations, see Zhu and Meng \cite{MR1775845}. See also e.g. Jiang \cite{MR1663321}, Shen and Jiang \cite{MR2240392} and Fu and Gao \cite{MR2563645} for related work on other infinite-dimensional Lie algebras, as well as Yang and Yu and Yao \cite{MR3057231}, and Fa and Han and Yue \cite{MR3514533}. Nonetheless, we will provide an explicit algebraic proof of $\mathrm{H}^1(\mathcal{W},\mathcal{W})=\{0\}$ in the appendix in order to illustrate our techniques.\\
The following result is a well-known fact:
\begin{equation*}
dim(\mathrm{H}^2(\mathcal{W},\mathbb{K}))=1\,,
\end{equation*}
which states that there is up to equivalence and rescaling exactly one non-trivial central extension of the Witt algebra, namely the Virasoro algebra.\\
Schlichenmaier \cite{MR3200354,MR3363999} and Fialowski \cite{MR2985241} showed that 
\begin{equation*}
\mathrm{H}^2(\mathcal{W},\mathcal{W})=\{0\}\,,
\end{equation*}
which implies infinitesimal and formal rigidity. 
However, contrary to finite-dimensional Lie algebras \cite{MR0171807,MR0207793,MR0240167,MR0178041}, the vanishing of $\mathrm{H}^2(\mathcal{W},\mathcal{W})$ does not imply rigidity with respect to other parameter spaces in the case of infinite-dimensional Lie algebras, see \cite{MR2030563,MR2363750,MR2183958,MR2509150,MR2269883}. \\
The main result of this article is the following:
\begin{equation*}
\mathrm{H}^3(\mathcal{W},\mathcal{W})=\{0\}\,,
\end{equation*}
see Theorem \ref{h3}.
Another result we obtain concerns the trivial module:
\begin{equation*}
dim (\mathrm{H}^3(\mathcal{W},\mathbb{K}))=1\,,
\end{equation*}
see Theorem \ref{h3K}.\\
For results concerning the Virasoro algebra $\mathcal{V}$, see Section \ref{Section6}. 



\section{Analysis of \texorpdfstring{$\mathrm{H}^3(\mathcal{W},\mathcal{W})$}{H3(W,W)}}\label{Section5}
In this section, we analyze the third cohomology group of the Witt algebra with values in the adjoint module.  The main aim of this section is to proof the following theorem:
\newtheorem{H3}{Theorem}[section]
\begin{H3} \label{h3}
The third cohomology of the Witt algebra $\mathcal{W}$ over a field $\mathbb{K}$ with char($\mathbb{K}$)$=0$ and values in the adjoint module vanishes, i.e.
\begin{equation*}
\mathrm{H}^3(\mathcal{W},\mathcal{W})=\{0\}\,.
\end{equation*}
\end{H3}
The proof follows in two steps, the first step concentrating on the non-zero degree part of the Witt algebra, the second step focusing on the degree zero part. \\
The coboundary condition for a 3-cocycle $\psi\in \mathrm{H}^3(\mathcal{W},\mathcal{W})$ is given by:
\begin{align*}
\psi(x_1,x_2,x_3)=(\delta_2\phi)(x_1,x_2,x_3)=&\phi\left(\left[x_1,x_2\right],x_3\right)+\phi\left(\left[x_2,x_3\right],x_1\right)+\phi\left(\left[x_3,x_1\right],x_2\right)\\
&-\left[x_1,\phi(x_2,x_3)\right]+\left[x_2,\phi(x_1,x_3)\right]-\left[x_3,\phi(x_1,x_2)\right]\,, 
\end{align*}
where $x_1,x_2,x_3\in\mathcal{W}$ and $\phi\in\mathrm{C}^2(\mathcal{W},\mathcal{W})$.\\
The cocycle condition for a 3-cocycle $\psi$ is given by:
\begin{align*}
&(\delta_3\psi)(x_1,x_2,x_3,x_4)\\
=&\psi\left(\left[x_1,x_2\right],x_3,x_4\right)-\psi\left(\left[x_1,x_3\right],x_2,x_4\right)+\psi\left(\left[x_1,x_4\right],x_2,x_3\right)\\
&+\psi\left(\left[x_2,x_3\right],x_1,x_4\right)-\psi\left(\left[x_2,x_4\right],x_1,x_3\right)
+\psi\left(\left[x_3,x_4\right],x_1,x_2\right)\\
&-\left[x_1,\psi(x_2,x_3,x_4)\right]+\left[x_2,\psi(x_1,x_3,x_4)\right]
-\left[x_3,\psi(x_1,x_2,x_4)\right]+\left[x_4,\psi(x_1,x_2,x_3)\right]=0\,.
\end{align*}
with $x_1,x_2,x_3,x_4\in\mathcal{W}$
\subsection{The non-zero degree part for the Witt algebra}
The proposition proved in this section shall be the following:
 \newtheorem{H3DegreeNonZero}{Proposition}[subsection]
\begin{H3DegreeNonZero} \label{H3DegreeNonZeroPart}
The following hold:
\begin{align*}
&\mathrm{H}^3_{(d)}(\mathcal{W},\mathcal{W})=\{0\}\ \text{ for }\ d\neq 0\ ,\\
&\mathrm{H}^3(\mathcal{W},\mathcal{W})=\mathrm{H}^3_{(0)}(\mathcal{W},\mathcal{W})\ .
\end{align*}
\end{H3DegreeNonZero}
\begin{proof}
Let $\psi\in\mathrm{H}^3_{(d\neq 0)}(\mathcal{W},\mathcal{W})$. \\
Let us perform a cohomological change $\psi'=\psi -\delta_2\phi$ with the following 2-cochain $\phi$:
\begin{equation*}
\phi(x_1,x_2)=-\frac{1}{d}\psi(x_1,x_2,e_0)\,,
\end{equation*}
which gives us, taking into account that $\phi(e_0,\cdot)=\phi(\cdot,e_0)=0$:
\begin{align*}
&\psi'(x_1,x_2,e_0)=\psi(x_1,x_2,e_0)-(\delta_2\phi)(x_1,x_2,e_0)\\
&=\psi(x_1,x_2,e_0)-\underbrace{\phi\left(\left[x_1,x_2\right],e_0\right)}_{=0}-\phi\left(\left[x_2,e_0\right],x_1\right)-\phi\left(\left[e_0,x_1\right],x_2\right)\\
&+[x_1,\underbrace{\phi(x_2,e_0)}_{=0}]-[x_2,\underbrace{\phi(x_1,e_0)}_{=0}]+\left[e_0,\phi(x_1,x_2)\right]\\
&=\psi(x_1,x_2,e_0)+\text{deg}(x_2)\underbrace{\phi(x_2,x_1)}_{=-\phi(x_1,x_2)}-\text{deg}(x_1)\phi(x_1,x_2)+\left(\text{deg}(x_1)+\text{deg}(x_2)+d\right)\phi(x_1,x_2)\\
&=-d\ \phi(x_1,x_2)+ d\ \phi(x_1,x_2)=0\,.
\end{align*}
We thus have \framebox{$\psi'(x_1,x_2,e_0)=0$}.\\
Next, let us write down the cocycle condition for $\psi'$ on the quadruplet $(x_1,x_2,x_3,e_0)$:
\begin{align*}
&(\delta_3\psi')(x_1,x_2,x_3,e_0)=0\\
\Leftrightarrow\  &\underbrace{\psi'\left(\left[x_1,x_2\right],x_3,e_0\right)}_{=0}-\underbrace{\psi'\left(\left[x_1,x_3\right],x_2,e_0\right)}_{=0}+\psi'\left(\left[x_1,e_0\right],x_2,x_3\right)\\
&+\underbrace{\psi'\left(\left[x_2,x_3\right],x_1,e_0\right)}_{=0}-\psi'\left(\left[x_2,e_0\right],x_1,x_3\right)
+\psi'\left(\left[x_3,e_0\right],x_1,x_2\right)\\
&-[x_1,\underbrace{\psi'(x_2,x_3,e_0)}_{=0}]+[x_2,\underbrace{\psi'(x_1,x_3,e_0)}_{=0}]
-[x_3,\underbrace{\psi'(x_1,x_2,e_0)}_{=0}]+\left[e_0,\psi'(x_1,x_2,x_3)\right]=0\\
\Leftrightarrow\  &-\text{deg}(x_1)\psi(x_1,x_2,x_3)+\text{deg}(x_2)\underbrace{\psi(x_2,x_1,x_3)}_{=-\psi(x_1,x_2,x_3)}-\text{deg}(x_3)\underbrace{\psi(x_3,x_1,x_2)}_{=\psi(x_1,x_2,x_3)}\\
&+\left(\text{deg}(x_1)+\text{deg}(x_2)+\text{deg}(x_3)+d\right)\psi(x_1,x_2,x_3)=0\\
\Leftrightarrow\ \ &d\ \psi(x_1,x_2,x_3)=0 \Leftrightarrow\  \psi(x_1,x_2,x_3)=0 \text{ as } d\neq 0\,.
\end{align*}
We conclude that the third cohomology of the Witt algebra reduces to the degree zero part, in agreement with the result of Fuks (\ref{Fuks}).
\end{proof}
\subsection{The degree zero part for the Witt algebra}
The proposition we shall prove in this section is the following:
 \newtheorem{H3DegreeZero}{Proposition}[subsection]
\begin{H3DegreeZero} \label{H3DegreeZeroPart}
The following holds:
\begin{equation*}
\mathrm{H}^3_{(0)}(\mathcal{W},\mathcal{W})=\{0\}\ .
\end{equation*}
\end{H3DegreeZero}
Clearly, Proposition \ref{H3DegreeZeroPart} together with Proposition \ref{H3DegreeNonZeroPart} shows Theorem \ref{h3}.
The proof of Proposition \ref{H3DegreeZeroPart} is accomplished in six steps and is similar to the proof performed for $\mathrm{H}^2_{(0)}(\mathcal{W},\mathcal{W})$ in \cite{MR3200354,MR3363999}. \\
Let $\psi$ be a degree zero 3-cocycle, i.e. we can write it as $\psi(e_i,e_j,e_k)=\psi_{i,j,k}e_{i+j+k}$ with suitable coefficients $\psi_{i,j,k}\in\mathbb{K}$. We say that $\psi_{\cdot,\cdot,\cdot}$ is of \textit{level} $l\in\mathbb{Z}$ if one of its indices is equal to $l$, i.e. $\psi_{\cdot,\cdot,\cdot}=\psi_{\cdot,\cdot,l}$ or some permutation thereof.\\ 
Consequently, five steps of the proof correspond to the analysis of the levels plus one, minus one, zero, plus two and minus two. The final step consists in the analysis of generic levels, which is obtained by induction. In each step, there are always three cases to consider depending on the signs of the indices. One of the three indices corresponds to the level and is fixed. In that case, the three cases to consider correspond to both remaining indices being negative, both being positive, or one being negative and one being positive. It does not matter which of the indices are chosen to be positive or negative, nor does it matter which one of the three indices is chosen to be fixed, because of the alternating property of the cochains.  
In the following, we provide a brief and superficial summary of the proof:
\begin{itemize}
	\item \textbf{Level plus one / minus one:} There is a cohomological change $\psi'=\psi-\delta_2\phi$, $\phi\in C^2_{(0)}(\mathcal{W},\mathcal{W})$ which allows to normalize to zero either the coefficients of level plus one or the coefficients of level minus one, depending on the signs of the two remaining indices. More precisely, we normalize $\psi'$ to $\psi_{i,j,-1}'=0$ if $i$ and $j$ are both positive and $\psi_{i,j,1}'=0$ else. \\
	The aim is to use the coboundary condition to produce recurrence relations which provide a consistent definition of $\phi$, i.e. of all the $\phi_{i,j}\ \forall\ i,j\in\mathbb{Z}$. Each degree of freedom given by some $\phi_{\cdot,\cdot}$ should be used to cancel some coefficient of the form $\psi_{\cdot,\cdot,1}$ or $\psi_{\cdot,\cdot,-1}$. In the case where both indices of $\phi_{i,j}$ have the same sign, the definition of the $\phi_{i,j}$'s can be obtained in a straightforward manner from the recurrence relations. In the case where the two indices are of opposite sign, poles occur in the recurrence relations, and the definition of the $\phi_{i,j}$'s has to be obtained in a somewhat roundabout manner. 
\item \textbf{Level zero:} For a cocycle $\psi$ normalized as described in the previous bullet point, the cocycle conditions imply $\psi_{i,j,0}=0\ \forall\ i,j\in\mathbb{Z}$.\\
 The cocycle conditions provide recurrence relations which allow to deduce the result immediately for $i$ and $j$ of the same sign. For $i$ and $j$ of different sign, the proof is an (almost) straightforward generalization of the proof of $\mathrm{H}^2(\mathcal{W},\mathcal{W})_{(0)}=\{0\}$ given in \cite{MR3200354,MR3363999}.
\item \textbf{Level minus one / plus one:} The cocycle conditions imply $\psi_{i,j,1}=0$ if $i$ and $j$ are both positive and $\psi_{i,j,-1}=0$ else. Together with the result of the first bullet point, we have $\psi_{i,j,1}=\psi_{i,j,-1}=0\ \forall\ i,j\in\mathbb{Z}$.\\
	This step is the simplest one of the entire proof. The cocycle conditions provide again recurrence relations which allow to deduce the results directly.
\item \textbf{Levels plus two and minus two / Generic Level $k$:} The cocycle conditions imply $\psi_{i,j,-2}=0$ and $\psi_{i,j,2}=0\ \forall\ i,j\in\mathbb{Z}$. Induction on $k$ subsequently implies $\psi_{i,j,k}=0\ \forall\ i,j,k\in\mathbb{Z}$.\\
For both indices $i$ and $j$ negative, the first step consists in proving that level minus two is zero, i.e. $\psi_{i,j,-2}=0$. Induction on the third index allows to conclude that the coefficients $\psi_{i,j,k}$ are zero for all negative indices $i,j,k\leq 0$. These results can be obtained directly from the recurrence relations given by the cocycle conditions.\\
In the case of one positive and one negative index, the first step consists in proving that both levels plus two and minus two are zero, $\psi_{i,j,2}=\psi_{i,j,-2}=0$. This has to be done by using induction on either $i$ or $j$ depending on the level under consideration. Note that in the proof of $\mathrm{H}^2(\mathcal{W},\mathcal{W})_{(0)}=\{0\}$ in \cite{MR3200354,MR3363999}, the vanishing of the levels plus two and minus two could be proved directly without using induction. Obviously, the number of times induction has to be used increases with the number of indices. Due to poles and zeros in the recurrence relations, the proof again follows a somewhat roundabout way.
The second and final step consists in using induction on the third index in order to prove $\psi_{i,j,k}=0$ for mixed indices, i.e. two indices positive and one index negative or two indices negative and one index positive. \\
The final case with both indices $i$ and $j$ positive starts with the proof that level plus two is zero, i.e. $\psi_{i,j,2}=0$. Induction on the third index allows to conclude that the coefficients $\psi_{i,j,k}$ are zero for all positive indices $i,j,k \geq  0$. These results follow directly from the recurrence relations.  
\end{itemize}
We now come to the detailed proof.
Let us write down the coboundary and cocycle conditions for later use. If $\phi$ is a degree zero 2-cochain, i.e. $\phi(e_i,e_j)=\phi_{i,j}e_{i+j}$, the coboundary condition for $\psi$ on the triplet $(e_i,e_j,e_k)$ becomes:
\begin{align*}
\psi_{i,j,k}=(\delta_2\phi)_{i,j,k}=&(j-i)\phi_{i+j,k}+(k-j)\phi_{k+j,i}+(i-k)\phi_{i+k,j}\\
&-(j+k-i)\phi_{j,k}+(i+k-j)\phi_{i,k}-(i+j-k)\phi_{i,j}\,.
\end{align*}
The cocycle condition for $\psi$ on the quadruplet $(e_i,e_j,e_k,e_l)$ becomes:
\begin{align*}
(\delta_3\psi)_{i,j,k,l}=&(j-i)\psi_{i+j,k,l}-(k-i)\psi_{i+k,j,l}+(l-i)\psi_{i+l,j,k}\\
&+(k-j)\psi_{k+j,i,l}-(l-j)\psi_{l+j,i,k}+(l-k)\psi_{l+k,i,j}\\
&-(j+k+l-i)\psi_{j,k,l}+(i+k+l-j)\psi_{i,k,l}\\
&-(i+j+l-k)\psi_{i,j,l}+(i+j+k-l)\psi_{i,j,k}=0\,.
\end{align*}
The first step of the proof is achieved with a cohomological change:
 \newtheorem{H3MinusOne}{Lemma}[subsection]
\begin{H3MinusOne} \label{Lemma1}
Every 3-cocycle $\psi$ of degree zero is cohomologous to a degree zero 3-cocycle $\psi'$ with: 
\begin{align}
\begin{array}[h]{lll}
	&\psi_{i,j,1}'=0& \forall\ i\leq 0,\ \forall\ j \in\mathbb{Z}\,, \\
	\text{  and  } &\psi_{i,j,-1}'=0& \forall\ i,j> 0\,,\\
\text{  and  }&\psi_{i,-1,2}'=0& \forall\ i\in\mathbb{Z}\,,  \\
\text{  and  }&\psi_{-4,2,-2}'=0\,. &
\end{array}\label{ResLemma1}
\end{align}
\end{H3MinusOne}
\begin{proof}
If $\phi$ is a 2-cochain, it can always be normalized to $\phi_{i,1}=0\ \forall\ i\in\mathbb{Z}$ and $\phi_{-1,2}=0$ with a cohomological change. The proof can be found in \cite{MR3200354,MR3363999} where it was performed in the context of $\phi$ being a 2-cocycle. However, as the cocycle condition is not used in the proof, the result is valid for any 2-cochain $\phi$. Hence, we will perform a cohomological change $\psi'=\psi-\delta_2\phi$ with $\phi$ normalized to \framebox{$\phi_{i,1}=0\ \forall\ i\in\mathbb{Z}$} and \framebox{$\phi_{-1,2}=0$}. This simplifies the notations considerably.\\
To increase the readability of the proof, we will separate the analysis depending on the signs of the indices $i,j$. Let us start with the case $i$ and $j$ both being negative.\\
\framebox{\textbf{Case 1:} $i,j\leq 0$}\\
Our aim is to show that we can find coefficients $\phi_{i,j}$ such that $\psi_{i,j,1}'=0$. 
Writing down the coboundary condition for $(i,j,1)$ and dropping the terms of the form $\phi_{.,1}$, we need:
\begin{equation*}
\psi_{i,j,1}=-(i+j-1)\phi_{i,j}+(i-1)\phi_{i+1,j}-(j-1)\phi_{j+1,i}\,. 
\end{equation*} 
This is the case if we define $\phi$:
\begin{equation*}
\phi_{i,j}:=\frac{i-1}{i+j-1}\phi_{i+1,j}-\frac{j-1}{i+j-1}\phi_{j+1,i}-\frac{\psi_{i,j,1}}{i+j-1}\,.
\end{equation*}
Starting with $i=0$, $j=-1$, $j$ decreasing and using $\phi_{.,1}=0$, this recurrence relation defines in a first step $\phi_{0,j}$ for $j\leq -1$. In a second step, $\phi_{-1,j}$ with $j\leq -2$ can be obtained, and so on for all $i\leq -2$ with $j<i$. It is sufficient to consider $j<i$ due to the alternating character of the cochains.
Thus, this recurrence relation defines $\phi_{i,j}$ for $i,j\leq 0$. It follows that we can perform a cohomological change such that $\psi_{i,j,1}'=0\ \forall\ i,j\leq 0$.\\
\framebox{\textbf{Case 2:} $i\leq 0$ and $j>0$}\\
We will start by proving that we can obtain $\psi_{i,2,-1}'=0\ \forall\ i \leq 0$ for a suitable choice of the coefficients $\phi_{i,j}$.\\
Let us consider the coboundary condition for $(-3,2,-1)$. Taking into account the normalization $\phi_{2,-1}=0$ we obtain:
\begin{equation*}
-2\phi_{-4,2}-6\phi_{-3,-1}=\psi_{-3,2,-1}\,.
\end{equation*}
The quantity $\phi_{-3,-1}$ has been defined in the previous case $i,j\leq 0$. Thus, we obtain a definition for $\phi_{-4,2}$. From there, we can obtain $\phi_{i,2}\ i\leq -5$ by using the coboundary condition for $(i,2,-1)$ and $\phi_{2,-1}=0$, which gives us:
\begin{equation}
\phi_{i-1,2}=\frac{3+i}{i+1}\phi_{i,2}+\frac{\psi_{i,2,-1}}{i+1}-\frac{i-3}{i+1}\phi_{i,-1}+\frac{i-2}{i+1}\phi_{i+2,-1}\,.\label{Phi2}
\end{equation}
The last two terms have been defined in the previous case $i,j\leq 0$. Thus, this defines $\phi_{i,2}\ i\leq -4$ such that we have $\psi'_{i,2,-1}=0\ \forall\ i\leq -3$. Next, let us consider the coboundary condition for $(-4,2,-2)$:
\begin{equation*}
\psi_{-4,2,-2}=-2\phi_{-6,2}-8\phi_{-4,-2}-4\phi_{0,-4}-4\phi_{2,-2}\,.
\end{equation*}
The coefficients $\phi_{-4,-2}$ and $\phi_{0,-4}$ have been defined in the previous case $i,j\leq 0$. The coefficient $\phi_{-6,2}$ has been defined in (\ref{Phi2}) for $i\leq -4$. Therefore, we obtain a definition for $\phi_{2,-2}$, which annihilates $\psi'_{-4,2,-2}$, \framebox{$\psi'_{-4,2,-2}=0$}.\\
As $\phi_{2,-2}$ is now defined, we can come back to Equation (\ref{Phi2}), insert $i=-2$ and obtain a definition for $\phi_{-3,2}$, annihilating $\psi'_{-2,2,-1}$. Since $\phi_{-1,2}=0$ due to our normalization, the only remaining $\phi_{i,2} \ i\leq 0$ to define is $\phi_{0,2}$.\\
Let us write down the coboundary condition for $(0,2,-1)$:
\begin{align*}
&-(3\phi_{0,2}+3\phi_{0,-1})=\psi_{0,2,-1}\\
\Leftrightarrow\  & \phi_{0,2}=-\phi_{0,-1}-\frac{1}{3}\psi_{0,2,-1}\,.
\end{align*}
This defines $\phi_{0,2}$ and consequently, $\psi_{0,2,-1}'=0$. Since $\psi_{-1,2,-1}'=0$ due to the alternating property, we obtain all in all that $\psi_{i,2,-1}'=0\ \forall\ i\leq 0$.\\
Next, let us prove that we can obtain $\psi_{i,j,1}'=0\ \forall\ i\leq 0\ \forall\ j>0$. It suffices to write down the coboundary condition for $(i,j,1)$ in the following way:
\begin{equation*}
\phi_{i,j+1}:=\frac{i+j-1}{j-1}\phi_{i,j}-\frac{i-1}{j-1}\phi_{i+1,j}+\frac{\psi_{i,j,1}}{j-1}\,.
\end{equation*}
Fixing $i=0$, and starting with $j=2$ (recall that $\phi_{i,1}=0$ and that we have just defined all $\phi_{i,2}\ i\leq 0$), $j$ increasing, we obtain $\phi_{0,j}\ \forall\ j>2$ and $\psi_{0,j,1}'=0\ \forall\ j\geq 2$. Similarly, fixing $i=-1$, and starting with $j=2$, $j$ increasing, we obtain $\phi_{-1,j}\ \forall\ j>2$ and $\psi_{-1,j,1}'=0\ \forall\ j\geq 2$. Continuing along the same lines, we obtain $\phi_{i,j}\ \forall\ i\leq 0,\ j>0$ and $\psi'_{i,j,1}=0\ \forall\ i\leq 0,\ j>0$. Together with the result $\psi_{i,j,1}'=0\ \forall\ i,j\leq 0$ obtained from the previous case with $i,j\leq 0$, we get \framebox{$\psi_{i,j,1}'=0\ \forall\ i\leq 0,\ \forall\ j \in\mathbb{Z}$} \\
\framebox{\textbf{Case 3:} $i> 0$ and $j>0$}\\  
Let us write down the coboundary condition for $(i,j,-1)$:
\begin{align*}
\psi_{i,j,-1}=&(i+1)\phi_{i-1,j}+(i-j-1)\phi_{i,-1}-(1+i+j)\phi_{i,j}\\
& +(j+1)\phi_{i,j-1}+(1+i-j)\phi_{j,-1}+(j-i)\phi_{i+j,-1}\,.
\end{align*}
From there, we can define $\phi$ via recurrence as follows:
\begin{align*}
\phi_{i,j}=&\frac{(i+1)}{(1+i+j)}\phi_{i-1,j}+\frac{(j+1)}{(1+i+j)}\phi_{i,j-1}-\frac{\psi_{i,j,-1}}{(1+i+j)}\\
&+\frac{(i-j-1)}{(1+i+j)}\phi_{i,-1}+\frac{(1+i-j)}{(1+i+j)}\phi_{j,-1}+\frac{(j-i)}{(1+i+j)}\phi_{i+j,-1}\,.
\end{align*}
Note that $\phi_{.,-1}$ have been defined in the previous case for $i\leq 0,\ j>0$. Starting with $i=2$, $j=3$ and $j$ increasing, we obtain in a first step $\phi_{2,j},\ \forall\ j\geq 3$ and $\psi'_{2,j,-1}=0\ \forall\ j\geq 3$. Next, fixing $i=3$, starting with $j=4$ and $j$ increasing, we obtain in a second step $\phi_{3,j},\ \forall\ j\geq 4$ and $\psi'_{3,j,-1}=0\ \forall\ j\geq 4$. Continuing similarly with $i$ increasing, we finally obtain all  $\phi_{i,j},\ \forall\ i,j>0$, and \framebox{$\psi'_{i,j,-1}=0\ \forall\ i,j>0$}. Note that we already have $\psi'_{1,j,-1}=0\ \forall\ j>0$ due to the previous case, which yielded $\psi_{i,j,1}'=0\ \forall\ i\leq 0,\ \forall\ j \in\mathbb{Z}$.
Combining the result $\psi'_{i,j,-1}=0\ \forall\ i,j>0$ with the result $\psi'_{i,2,-1}=0\ \forall\ i\leq 0$ from the previous case $i\leq 0,\ j>0$, we also obtain \framebox{$\psi'_{i,2,-1}=0\ \forall\ i\in\mathbb{Z}$}.
\end{proof}
\newtheorem{H3Zero}[H3MinusOne]{Lemma}
\begin{H3Zero} \label{Lemma2}
Let $\psi$ be a degree zero 3-cocycle such that: 
\begin{align*}
\begin{array}[h]{lll}
&\psi_{i,j,1}=0& \forall\ i\leq 0,\ \forall\ j \in\mathbb{Z} \,,\\
\text{  and  }&\psi_{i,j,-1}=0& \forall\ i,j> 0\,,\\
\text{  and  }&\psi_{i,-1,2}=0& \forall\ i\in\mathbb{Z}\,,
\end{array}
\end{align*}
then  
\begin{equation}
\psi_{i,j,0}=0\quad \forall\ i,j\in\mathbb{Z}\,.\label{ResLemma2}
\end{equation}
\end{H3Zero}
\begin{proof}
Again, we split the proof into the three cases depending on the signs of $i$ and $j$.\\
\framebox{\textbf{Case 1:} $i,j\leq 0$}\\
Let us write down the cocycle condition for $(i,j,0,1)$, neglecting the terms of the form $\psi_{i,j,1}\ i,j\leq 0$:
\begin{equation*}
(i+j-1)\psi_{i,j,0}-(i-1)\psi_{i+1,j,0}+(j-1)\psi_{j+1,i,0}=0\,.
\end{equation*}
We can define the following recurrence relation for $i$ and $j$ decreasing:
\begin{equation*}
\psi_{i,j,0}= \frac{(i-1)}{(i+j-1)}\psi_{i+1,j,0}-\frac{(j-1)}{(i+j-1)}\psi_{j+1,i,0}\,.
\end{equation*}
Fixing $i=-1$, starting with $j=-2$ and $j$ decreasing, we obtain $\psi_{-1,j,0}=0\ \forall\ j\leq -2$. Repeating the same procedure with decreasing values for $i$ and $j<i$, we obtain $\psi_{i,j,0}=0\ \forall\ i,j\leq 0$. \\
\framebox{\textbf{Case 2:} $i\leq 0$, $j>0$}\\
Let us write down the cocycle condition for $(i,2,0,-1)$:
\begin{align*}
&-\cancel{\psi_{-1,i,2}}+3 \psi_{1,i,0}+(-1+i) \psi_{2,0,-1}-2 \cancel{\psi_{2,i,-1}}-(1+i) \psi_{-1+i,2,0}\\
&+(-3+i) \psi_{i,0,-1}-\cancel{\psi_{i,2,-1}}+(3+i) \psi_{i,2,0}-(-2+i) \psi_{2+i,0,-1}=0\,.
\end{align*}
The slashed terms cancel each other, although they are zero anyway as we have $\psi_{i,2,-1}=0\ \forall\ i\in\mathbb{Z}$. The term $\psi_{1,i,0}$ is zero as we have $\psi_{i,j,1}=0\ \forall\ i,j\leq 0$. The term $\psi_{2,0,-1}$ is zero due to $\psi_{i,2,-1}=0\ \forall\ i\in\mathbb{Z}$. The terms $\psi_{i,0,-1}$ and $\psi_{2+i,0,-1}$ (for $i\leq -2$) are  zero because of the previous case, $\psi_{i,j,0}=0\ \forall\ i,j\leq 0$. 
Therefore, we are left with:
\begin{equation}
\psi_{i-1,2,0}=\frac{i+3}{i+1}\psi_{i,2,0}\,.\label{TwoZero}
\end{equation}
Putting $i=-3$ in the equation above, this recurrence relation implies $\psi_{-4,2,0}=0$ and by recursion $\psi_{i,2,0}=0\ \forall i\leq -4$. Next, consider the cocycle condition for $(i,2,-2,0)$:
\begin{align*}
&2 \cancel{\psi_{-2,i,2}}+i \psi_{2,-2,0}+2 \cancel{\psi_{2,i,-2}}+(2+i) \psi_{-2+i,2,0}\\
&+(-4+i) \psi_{i,-2,0}-(4+i) \psi_{i,2,0}-(-2+i) \psi_{2+i,-2,0}=0\,.
\end{align*}
The slashed terms cancel each other, the terms $\psi_{i,-2,0}$ and $\psi_{2+i,-2,0}$ (for $i\leq -2$) are zero because of $\psi_{i,j,0}=0\ \forall\ i,j\leq 0$. As we have $\psi_{i,2,0}=0\ \forall i\leq -4$, we can put for example $i=-4$ in the equation above and obtain $\psi_{2,-2,0}=0$. Inserting this value in Equation (\ref{TwoZero}) with $i=-2$, we obtain $\psi_{-3,2,0}=0$. Recall that we also have $\psi_{-1,2,0}=0$ due to $\psi_{i,-1,2}=0\ \forall\ i\in\mathbb{Z}$.
All in all, we have $\psi_{i,2,0}\ \forall\ i\leq 0$.\\
This result is needed to write down a well-defined recurrence relation. Writing down the cocycle condition for $(i,j,0,1)$ and neglecting the terms of the form $\psi_{i,j,1}$ with $i,j\leq 0$ and $i\leq 0\ ,j>0$, we obtain the following recurrence relation:
\begin{equation*}
\psi_{i,j+1,0}=\frac{i+j-1}{j-1}\psi_{i,j,0}-\frac{i-1}{j-1}\psi_{i+1,j,0}\,.
\end{equation*} 
Fixing $i=-1$, one starts with $j=2$ (since we already have $\psi_{i,j,1}=0\ i,j\leq 0$ and $\psi_{i,2,0}=0\ i\leq 0$), which gives, with increasing $j$, $\psi_{-1,j,0}=0\ j\geq 3$. Continuing with fixing $i=-2$, starting again with $j=2$ and increasing $j$, we obtain $\psi_{-2,j,0}=0\ j\geq 3$. Doing this for all $i\leq 0$, we finally obtain $\psi_{i,j,0}=0\ \forall\ i\leq 0,\ j>0$.\\
 \framebox{\textbf{Case 3:} $i,j>0$}\\
Writing down the cocycle condition for $(i,j,0,-1)$, we obtain:
\begin{align*}
&-\cancel{\psi_{-1,i,j}}-(1+i) \psi_{-1+i,j,0}+(-1+i-j) \psi_{i,0,-1}+i \cancel{\psi_{i,j,-1}}-\cancel{(-1+i+j)\psi_{i,j,-1}}\\
&+(1+i+j) \psi_{i,j,0}+(1+j) \psi_{-1+j,i,0}+(1+i-j) \psi_{j,0,-1}-j \cancel{\psi_{j,i,-1}}+(-i+j) \psi_{i+j,0,-1}=0\,.
\end{align*}
The slashed terms cancel each other, though they are zero anyway due to $\psi_{i,j,-1}=0\ i,j>0$. The terms $\psi_{i,0,-1}$, $\psi_{j,0,-1}$ and $\psi_{i+j,0,-1}$ are zero due to the previous case, $\psi_{i,j,0}=0\ i\leq 0,\ j>0$. Thus, we obtain the following recurrence relation:
\begin{equation*}
\psi_{i,j,0}=\frac{(1+i)}{(1+i+j)} \psi_{-1+i,j,0}-\frac{(1+j)}{(1+i+j)} \psi_{-1+j,i,0}\,.
\end{equation*} 
Fixing $i=1$, starting with $j=2$, $j$ increasing, we obtain $\psi_{1,j,0}=0\ j\geq 2$. Fixing $i=2$, starting with $j=3$, we get $\psi_{2,j,0}=0\ j\geq 3$. Continuing with increasing $i$ and keeping $j>i$ due to skew-symmetry, we finally obtain $\psi_{i,j,0}=0\ \forall\ i,j>0$.\\
Taking all three cases together, we obtain the announced result, $\psi_{i,j,0}=0\ \forall\ i,j\in\mathbb{Z}$. 
\end{proof}
\newtheorem{H3PlusOne}[H3MinusOne]{Lemma}
\begin{H3PlusOne} \label{Lemma3}
Let $\psi$ be a degree zero 3-cocycle such that: 
\begin{align*}
\begin{array}[h]{lll}
&\psi_{i,j,1}=0& \forall\ i\leq 0,\ \forall\ j \in\mathbb{Z}\,, \\ 
\text{  and  }&\psi_{i,j,-1}=0& \forall\ i,j> 0\,,\\
\text{  and  }&\psi_{i,-1,2}=0& \forall\ i\in\mathbb{Z} \,,\\
\text{  and  }&\psi_{i,j,0}=0& \forall\ i,j\in\mathbb{Z}\,,
\end{array}
\end{align*}
then  
\begin{equation*}
\psi_{i,j,1}=\psi_{i,j,-1}=0\quad \forall\ i,j\in\mathbb{Z}\,.
\end{equation*}
\end{H3PlusOne}
\begin{proof}
Again, the proof is split into the three cases depending on the signs of $i,j$.\\
\framebox{\textbf{Case 1:} $i,j\leq 0$}\\
Writing down the cocycle condition for $(i,j,1,-1)$ and neglecting $\psi_{i,j,1}\ i,j\leq 0$ as well as $\psi_{i,j,0}\ i,j\leq 0$, we obtain:
\begin{align*}
&-(-2+i+j) \psi_{i,j,-1}+(-1+i) \psi_{1+i,j,-1}-(-1+j) \psi_{1+j,i,-1}=0\\
\Leftrightarrow\  & \psi_{i,j,-1}=\frac{(-1+i)}{(-2+i+j)} \psi_{1+i,j,-1}-\frac{(-1+j)}{(-2+i+j)} \psi_{1+j,i,-1}\,.
\end{align*}
Fixing $i=-2$ (since level zero $\psi_{0,j,-1}\ j\leq 0$ is already done and $\psi_{-1,j,-1}=0$), starting with $j=-3$ and $j$ decreasing, we obtain $\psi_{-2,j,-1}=0\ j\leq -3$. Fixing $i=-3$, starting with $j=-4$ and $j$ decreasing,  we get $\psi_{-3,j,-1}=0\ j\leq -4$. Continuing along the same lines, we obtain $\psi_{i,j,-1}=0\ i,j\leq 0$.\\
 \framebox{\textbf{Case 2:} $i\leq 0,\ j>0$}\\
Writing down the cocycle condition for $(i,j,1,-1)$ and neglecting $\psi_{i,j,1}\ i\leq 0,\ j>0$, $\psi_{i,j,1}\ i,j\leq 0$ as well as $\psi_{i,j,0}\ i\leq 0,\ j>0$, we obtain:
\begin{align*}
&-(-2 + i + j) \psi_{i, j, -1}+(-1+i) \psi_{1+i,j,-1}-(-1 + j) \psi_{1 + j, i, -1}=0\\
&\Leftrightarrow\   \psi_{ i,1 + j, -1}=\frac{(-2 + i + j)}{(-1 + j)} \psi_{i, j, -1}-\frac{(-1+i)}{(-1 + j)} \psi_{1+i,j,-1}\,.
\end{align*}
Fixing $i=-2$ (since level zero $\psi_{0,j,-1}\ j>0$ is already done and $\psi_{-1,j,-1}=0$) and starting with $j=2$ (since $\psi_{i,2,-1}=0\ i\leq 0$), increasing $j$, we obtain $\psi_{-2,j,-1}=0\ j\geq 3$. Fixing $i=-3$, starting again with $j=2$, $j$ increasing, we get $\psi_{-3,j,-1}=0\ j\geq 3$. Continuing with $i$ decreasing, we get $\psi_{i,j,-1}=0\ i\leq 0,\ j>0$.\\
\framebox{\textbf{Case 3:} $i,j>0$}\\
Writing again down the cocycle condition for $(i,j,1,-1)$, this time neglecting the terms $\psi_{i,j,-1}\ i,j>0$ and $\psi_{i,j,0}\ i,j>0$, we obtain:
\begin{align*}
&-(1 + i) \psi_{-1 + i, j, 1}+(2 + i + j) \psi_{i, j, 1}+(1 + j) \psi_{-1 + j, i, 1}=0\\
&\Leftrightarrow\  \psi_{i, j, 1}= \frac{(1 + i)}{(2 + i + j)} \psi_{-1 + i, j, 1}-\frac{(1 + j)}{(2 + i + j)} \psi_{-1 + j, i, 1}\,.
\end{align*}
Fixing $i=2$ (since $\psi_{1,j,1}=0$) and starting with $j=3$, increasing $j$, we obtain $\psi_{2,j,1}=0\ j\geq 3$. Increasing $i$ and keeping $j>i$ we finally obtain $\psi_{i,j,1}=0\ \forall\ i,j>0$.\\
Taking all three cases together, we have proven that $\psi_{i,j,1}=0\ \forall\ i,j\in\mathbb{Z}\text{  and  }\psi_{i,j,-1}=0\ \forall\ i,j\in\mathbb{Z}$.
\end{proof}
\newtheorem{H3PlusTwo}[H3MinusOne]{Lemma}
\begin{H3PlusTwo} \label{Lemma4}
Let $\psi$ be a degree zero 3-cocycle such that: 
\begin{equation*}
\psi_{i,j,1}=\psi_{i,j,-1}=\psi_{i,j,0}=0\quad \forall\ i,j\in\mathbb{Z}\quad \text{  and  }\quad \psi_{-4,2,-2}=0\,,
\end{equation*}
then  
\begin{equation*}
\psi_{i,j,k}=0\quad \forall\ i,j,k\in\mathbb{Z}\,.
\end{equation*}
\end{H3PlusTwo}
\begin{proof}
Again, the proof is split in the three cases depending on the signs of $i$ and $j$.\\
\framebox{\textbf{Case 1:} $i,j\leq 0$}\\
In a first step, we shall prove the following statement: \framebox{$\psi_{i,j,-2}=0\ \forall\ i,j\leq 0$}.\\
Writing down the cocycle condition for $(i,j,-2,1)$ and neglecting the terms of level one $\psi_{i,j,1}$ and level minus one $\psi_{i,j,-1}$, we obtain the following:
\begin{align*}
&(-3+i+j) \psi_{i,j,-2}-(-1+i) \psi_{1+i,j,-2}+(-1+j) \psi_{1+j,i,-2}=0\\
&\Leftrightarrow\  \psi_{i,j,-2}=\frac{(-1+i)}{(-3+i+j)} \psi_{1+i,j,-2}-\frac{(-1+j)}{(-3+i+j)} \psi_{1+j,i,-2}\,.
\end{align*}
Fixing $i=-3$ (since the levels zero $\psi_{i,j,0}$ and minus one $\psi_{i,j,-1}$ are already done and $\psi_{-2,j,-2}=0$), starting with $j=-4$ and decreasing $j$, we obtain $\psi_{-3,j,-2}=0\ j\leq -4$. Continuing along the same lines with decreasing $i$ and keeping $j<i$, we obtain $\psi_{i,j,-2}=0\ \forall\ i,j\leq 0$ as a first step. \\
In a second step, we shall prove \framebox{$\psi_{i,j,k}=0\ \forall\ i,j,k\leq 0$}. This can be done by induction. We know the result is true for $k=0,-1,-2$. Hence, we will assume it is true for some $k\leq -2$ and check whether it remains true for $k-1$. The cocycle condition for $(i,j,k,-1)$ is given by, after omitting terms of level minus one $\psi_{i,j,-1}$:
\begin{align*}
&-(1+i) \psi_{-1+i,j,k}+(1+i+j+k) \psi_{i,j,k}+(1+j) \psi_{-1+j,i,k}-(1+k) \psi_{-1+k,i,j}=0\\
& \Leftrightarrow\  -(1+k) \psi_{-1+k,i,j}=0\Leftrightarrow\  \psi_{-1+k,i,j}=0\text{ as }k\leq -2\,.
\end{align*}
The terms $\psi_{-1+i,j,k}$, $\psi_{i,j,k}$ and $\psi_{-1+j,i,k}$ are zero since they are of level $k$ and thus zero by induction hypothesis. It follows $\psi_{i,j,k}=0\ \forall\ i,j,k\leq 0$.\\
\framebox{\textbf{Case 2:} $i\leq 0,\ j>0$}\\
In a first step, we shall prove the following statements: \framebox{$\psi_{i,j,-2}=0\ \forall\ i\leq 0,\ j>0$} and \framebox{$\psi_{i,j,2}=0\ \forall\ i\leq 0,\ j>0$}.\\
The cocycle condition for $(-3,2,-2,-1)$ reads, after dropping the terms of level one, minus one and zero:
\begin{equation*}
2 \psi_{-4,2,-2}-2 \psi_{-3, 2, -2}=0\,.
\end{equation*}
Since we have $\psi_{-4,2,-2}=0$, the equation above implies $\psi_{-3, 2, -2}=0$. Next, let us write down the cocycle condition for $(-3,j,-2,1)$, which gives after dropping terms of level one and of level minus one:
\begin{align*}
&(-6 + j) \psi_{-3, j, -2}+(-1 + j) \psi_{1 + j, -3, -2}=0\\
&\Leftrightarrow\  \psi_{1 + j, -3, -2}=\frac{(-6 + j)}{(-1 + j)} \psi_{j, -3, -2}\,.
\end{align*}
Starting with $j=2$, we obtain $\psi_{j,-3,-2}=0\ \forall\ j\geq 3$ since the starting point is zero: $\psi_{2, -3, -2}=0$. Adding the level one, we obtain $\psi_{j,-3,-2}=0\ \forall\ j> 0$.\\
Next, let us write down the cocycle condition for $(i,3,2,-1)$ after dropping terms of level one and level minus one:
\begin{align}
&-(1 + i) \psi_{-1 + i, 3, 2}+(6 + i) \psi_{i, 3, 2}=0 \nonumber \\
&\Leftrightarrow\  \psi_{-1 + i, 3, 2}=\frac{(6 + i)}{(1 + i)} \psi_{i, 3, 2}\,.\label{level32}
\end{align}
This gives us for $i=-2$: $\psi_{-3, 3, 2}=-4 \psi_{-2, 3, 2}$\\
For $i=-3$: $\psi_{-4, 3, 2}=\frac{3}{-2} \psi_{-3, 3, 2}=6\psi_{-2, 3, 2}$\\
For $i=-4$: $\psi_{-5, 3, 2}=\frac{2}{-3} \psi_{-4, 3, 2}=-4\psi_{-2, 3, 2} $\\
Now, let us write down the cocycle condition for $(-3,3,2,-2)$ and drop the terms of level one, level minus one and level zero, as well as the terms of the form $\psi_{j,-3,-2}\ j>0$:
\begin{align*}
&\psi_{-5, 3, 2}+4 \psi_{-3, 3, 2}-6 \psi_{3, 2, -2}=0\\
&\Leftrightarrow\  -4\psi_{-2, 3, 2}-16 \psi_{-2, 3, 2}-6 \psi_{-2, 3, 2}=0\\
&\Leftrightarrow\  \psi_{-2, 3, 2}=0\,.
\end{align*}
To obtain the second line, the first two terms were simply replaced by their expressions computed above. Putting $i=-2$ and $\psi_{-2, 3, 2}=0$ in the recurrence relation (\ref{level32}), we obtain $\psi_{i, 3, 2}=0\ \forall\ i\leq -3$. Together with the levels minus one and zero, we obtain $\psi_{i, 3, 2}=0\ \forall\ i\leq 0$.\\
To prove $\psi_{i,j,-2}=0\ \forall\ i\leq 0,\ j>0$, we will use induction on $i$. Indeed, we have proven that $\psi_{i,j,-2}=0\ \forall\  j>0$, for $i=0,-1,-2,-3$ (recall that we have $\psi_{j,-3,-2}=0\ j>0$). Suppose the statement holds true 
down to $i+1$, $i \leq -4$, and let us see what happens for $i$. 
The cocycle condition for $(i,j,-2,1)$ gives, after dropping terms of level one and level minus one:
\begin{align}
&(-3+i+j) \psi_{i,j,-2}-(-1+i) \underbrace{\psi_{1+i,j,-2}}_{=0}+(-1+j) \psi_{1+j,i,-2}=0\nonumber \\
&\Leftrightarrow\  \psi_{i,j+1,-2}=\frac{(-3+i+j)}{(-1+j)} \psi_{i,j,-2}\,. \label{InductionEnd}
\end{align}
The term in the middle is zero due to the induction hypothesis.\\
This gives, for $j=2$: $\psi_{i,3,-2}=(-1+i) \psi_{i,2,-2}$.\\
For $j=3$: $\psi_{i,4,-2}=\frac{i}{2} \psi_{i,3,-2}=\frac{(-1+i)(i)}{2}\psi_{i,2,-2}$.\\
For $j=4$: $\psi_{i,5,-2}=\frac{(1+i)}{3} \psi_{i,4,-2}=\frac{(-1+i)(i)(1+i)}{6}\psi_{i,2,-2}$.\\
Next, we will insert these values into the cocycle condition for $(i,3,-2,2)$, after dropping terms of level zero and level one:
\begin{align*}
&(-3+i) \psi_{3,-2,2}+\psi_{5,i,-2}+(2+i) \psi_{-2+i,3,2}+(-3+i) \psi_{i,-2,2}\\ 
&+(-1+i)\psi_{i,3,-2}-(7+i) \psi_{i,3,2}-(-2+i) \psi_{2+i,3,-2}-(-3+i) \psi_{3+i,-2,2}=0\\
& \Leftrightarrow\  -\frac{(-1+i)(i)(1+i)}{6}\psi_{i,2,-2}-(-3+i) \psi_{i,2,-2}+(-1+i)(-1+i) \psi_{i,2,-2}=0\\
& \Leftrightarrow\  (i-3)(i^2-3i+8) \psi_{i,2,-2}=0\,.
\end{align*}
The terms $\psi_{3,-2,2}$, $\psi_{-2+i,3,2}$ and $\psi_{i,3,2}$ are zero due to what was proved before, $\psi_{i,3,2}=0\ \forall\ i\leq 0$. The terms $\psi_{2+i,3,-2}$ and $\psi_{3+i,-2,2}$ are zero as a consequence of the induction hypothesis. In the last line, we have $(i-3)\neq 0$, since $i\leq -4$, and also $(i^2-3i+8)\neq 0$ since its discriminant is negative. It follows $\psi_{i,2,-2}=0$. Reinserting this into (\ref{InductionEnd}) and taking into account that level one is zero, we obtain that the induction holds true for $i$, and thus: $\psi_{i,j,-2}=0\ \forall\ i\leq 0, j>0$.\\
Next, we proceed similarly, but with induction on $j$, to prove $\psi_{i,j,2}=0\ \forall\ i\leq 0, j>0$. We already know that the statement holds true for $j=1,2,3$. Let us suppose it is true up to $j-1$, $j\geq 4$, and show that it remains true for $j$. Let us write down the cocycle condition for $(i,j,2,-1)$, after dropping terms of level one and level minus one:
\begin{align}
&-(1+i) \psi_{-1+i,j,2}+(3+i+j) \psi_{i,j,2}+(1+j) \underbrace{\psi_{-1+j,i,2}}_{=0}=0\nonumber \\
&\Leftrightarrow\   \psi_{-1+i,j,2}=\frac{(3+i+j)}{(1+i)} \psi_{i,j,2}\,. \label{InductionEnd2}
\end{align}
The third term is zero due to the induction hypothesis.
From the recurrence relation above, we obtain for $i=-2$: $\psi_{-3,j,2}=-(1+j) \psi_{-2,j,2}$.\\
For $i=-3$: $\psi_{-4,j,2}=\frac{j}{(-2)} \psi_{-3,j,2}=\frac{j(1+j)}{2} \psi_{-2,j,2}$\\
For $i=-4$: $\psi_{-5,j,2}=\frac{(-1+j)}{(-3)} \psi_{-4,j,2}=-\frac{(-1+j)j(1+j)}{6}\psi_{-2,j,2}$.\\
Next, we insert these values into the cocycle condition for $(-3,j,2,-2)$ after dropping terms of level zero and level minus one:
\begin{align*}
&\psi_{-5,j,2}-(3+j) \psi_{-3,2,-2}-(-7+j) \psi_{-3,j,-2}+(1+j) \psi_{-3,j,2}\\
&+(3+j) \psi_{-3+j,2,-2}+(2+j) \psi_{-2+j,-3,2}-(3+j) \psi_{j,2,-2}-(-2+j) \psi_{2+j,-3,-2}=0\\
&\Leftrightarrow\  -\frac{(-1+j)j(1+j)}{6}\psi_{-2,j,2}-(1+j)(1+j) \psi_{-2,j,2}-(3+j) \psi_{-2,j,2}=0\\
&\Leftrightarrow\  (j+3)(8+3j+j^2)\psi_{-2,j,2}=0\,.
\end{align*}
The terms $\psi_{-3,2,-2}$, $\psi_{-3,j,-2}$ and $\psi_{2+j,-3,-2}$ are zero due to what was shown before in this proof, $\psi_{j,-3,-2}=0\ \forall\ j>0$. The terms $\psi_{-3+j,2,-2}$ and $\psi_{-2+j,-3,2}$ are zero due to the induction hypothesis. In the last line, we have $(j+3)\neq 0$, since $j\geq 4$, and also $(j^2+3j+8)\neq 0$ since its discriminant is negative. It follows $\psi_{-2,j,2}=0$. Reinserting this into (\ref{InductionEnd2}) and taking into account that level minus one is zero, we obtain that the induction holds true for $j$, and thus: $\psi_{i,j,2}=0\ \forall\ i\leq 0, j>0$. \\
Now that the terms of level two and of level minus two are zero for the case $i\leq 0, j>0$, we can use induction on $k$ to first prove \framebox{$\psi_{i,j,k}=0,\ \forall\ i\leq 0, j>0, k\geq 0$} and then \\ \framebox{$\psi_{i,j,k}=0,\ \forall\ i\leq 0, j>0, k\leq 0$}. \\
The result is true for $k=0,1,2$. Let us assume the result is true for $k$, $k\geq 2$ and show that it remains true for $k+1$. The cocycle condition for $(i,j,k,1)$ gives, after dropping terms of level one,
\begin{align*}
&(-1+i+j+k) \psi_{i,j,k}-(-1+i) \psi_{1+i,j,k}+(-1+j) \psi_{1+j,i,k}-(-1+k) \psi_{1+k,i,j}=0\\
&\Leftrightarrow\  (-1+k) \psi_{1+k,i,j}=0\Leftrightarrow\  \psi_{1+k,i,j}=0\text{ as }k\geq 2\,.
\end{align*}
The terms $\psi_{i,j,k}$, $\psi_{1+i,j,k}$ and $\psi_{1+j,i,k}$ are zero because of the induction hypothesis (if $1+i=1$, the term $\psi_{1+i,j,k}$ is still zero because the level plus one is zero for all $j,k\in\mathbb{Z}$). It follows that the result holds true for $k+1$. \\
All the same, the result is true for $k=0,-1,-2$. Let us assume it is true for $k$, $k\leq -2$, and show that it holds true for $k-1$. The cocycle condition for $(i,j,k,-1)$ yields, after dropping terms of level minus one,
\begin{align*}
&-(1+i) \psi_{-1+i,j,k}+(1+i+j+k) \psi_{i,j,k}+(1+j) \psi_{-1+j,i,k}-(1+k) \psi_{-1+k,i,j}=0\\
&\Leftrightarrow\  (1+k) \psi_{-1+k,i,j}=0\Leftrightarrow\  \psi_{-1+k,i,j}=0\text{ as }k\leq -2\,.
\end{align*}
The terms $\psi_{-1+i,j,k}$, $\psi_{i,j,k}$ and $\psi_{-1+j,i,k}$ are zero because of the induction hypothesis (if $-1+j=0$, the term $\psi_{-1+j,i,k}$ is still zero because the level zero vanishes for all $i,k\in\mathbb{Z}$). It follows that the result holds true for $k-1$. Thus, we have obtained the desired result for the case $i\leq 0, j>0$.\\
\framebox{\textbf{Case 3:} $i>0, j>0$}\\
In a first step, we shall prove the following statement: \framebox{$\psi_{i,j,2}=0\ \forall\ i,j>0$}. The cocycle condition for $(i,j,2,-1)$ yields, after dropping the terms of level one and of level minus one:
\begin{align*}
&-(1+i) \psi_{-1+i,j,2}+(3+i+j) \psi_{i,j,2}+(1+j) \psi_{-1+j,i,2}=0\\
&\Leftrightarrow\  \psi_{i,j,2}=\frac{(1+i)}{(3+i+j)} \psi_{-1+i,j,2}+\frac{(1+j)}{(3+i+j)} \psi_{i,-1+j,2}\,.
\end{align*} 
Fixing $i=3$ and starting with $j=4$, $j$ ascending, we obtain $\psi_{3,j,2}=0\ \forall j\geq 4$. Fixing $i=4$ and starting with $j=5$, $j$ ascending, we get $\psi_{4,j,2}=0\ \forall j\geq 5$. Continuing with ascending $i$, keeping $j>i$, we finally obtain $\psi_{i,j,2}=0\ \forall\ i,j>0$. \\
Finally, we want to prove \framebox{$\psi_{i,j,k}=0\ \forall\ i,j>0, k\geq 0$}. This can be done with induction on $k$. Indeed, the result is true for level zero, level one and level two, i.e. $k=0,1,2$. Thus, let us assume the result is true for $k$, $k\geq 2$ and show that it holds true for $k+1$. The cocycle condition for $(i,j,k,1)$ gives, after dropping the terms of level one:
\begin{align*}
&(-1+i+j+k)\psi_{i,j,k}-(-1+i) \psi_{1+i,j,k}+(-1+j) \psi_{1+j,i,k}-(-1+k) \psi_{1+k,i,j}=0\\
& \Leftrightarrow\  (-1+k) \psi_{1+k,i,j}=0\Leftrightarrow\   \psi_{1+k,i,j}=0\text{ as }k\geq 2\,.
\end{align*}
The terms $\psi_{i,j,k}$, $\psi_{1+i,j,k}$ and $\psi_{1+j,i,k}$ are zero because of the induction hypothesis. It follows that the statement holds true for $k+1$.\\
Taking all three cases together, we find the announced result \framebox{$\psi_{i,j,k}=0\ \forall\ i,j,k\in\mathbb{Z}$}
\end{proof}
\textbf{Proof of the Proposition \ref{H3DegreeZeroPart}}\\
\begin{proof}
Let us collect the statements of the four lemmata. Let $\psi$ be a degree-zero 3-cocycle of $\mathcal{W}$ with values in $\mathcal{W}$. By Lemma \ref{Lemma1} we can perform a cohomological change such that we obtain a cohomologous degree-zero 3-cocycle with coefficients fulfilling (\ref{ResLemma1}). Hence, the assumptions of Lemma \ref{Lemma2} are satisfied and we obtain (\ref{ResLemma2}). Together with Lemma \ref{Lemma3}, the assumptions of Lemma \ref{Lemma4} are fulfilled and Lemma \ref{Lemma4} shows \framebox{$\psi_{i,j,k}=0\ \forall\ i,j,k\in\mathbb{Z}$}, which proves the Proposition \ref{H3DegreeZeroPart}.
\end{proof}
\section{The cohomology of the Virasoro algebra}\label{Section6}
In this section, we exhibit results concerning the cohomology of the Virasoro algebra, including known results as well as new results obtained in this article. The results are summarized in the first subsection, a sketch of the proofs of the new results are given in the second subsection. Details will be provided in forthcoming work \cite{Ecker}. 
\subsection{Results for the cohomology of the Virasoro algebra}\label{ResVirasoro}
From Section \ref{Zeroth cohomology} we immediately deduce the following results for the zeroth cohomology group of the Virasoro algebra:
\begin{equation*}
\mathrm{H}^0(\mathcal{V},\mathbb{K})=\mathbb{K}\qquad\text{ and }\qquad\mathrm{H}^0(\mathcal{V},\mathcal{V})=\mathbb{K}\ t\,, 
\end{equation*}
where $t$ is the central element. Concerning the first cohomology group, we obtain from Section \ref{First cohomology} the following results:
\begin{equation*}
\mathrm{H}^1(\mathcal{V},\mathbb{K})=\{0\} \qquad\text{ and }\qquad\mathrm{H}^1(\mathcal{V},\mathcal{V})=\{0\}\,. 
\end{equation*}
The first equation is obtained from the fact that the Virasoro algebra is a perfect Lie algebra. The second equation states that all the derivations for $\mathcal{V}$ are inner. This might be known, but nonetheless we will provide an explicit proof of $\mathrm{H}^1(\mathcal{V},\mathcal{V})=\{0\}$ in the appendix. \\
The results for the second cohomology were shown by Schlichenmaier \cite{MR3200354}:
\begin{equation*}
\mathrm{H}^2(\mathcal{V},\mathbb{K})=\{0\} \qquad\text{ and }\qquad\mathrm{H}^2(\mathcal{V},\mathcal{V})=\{0\}\,. 
\end{equation*}
Thus, for the second cohomology we have $\mathrm{H}^2(\mathcal{W},\mathcal{W})=\mathrm{H}^2(\mathcal{V},\mathcal{V})=\{0\}$. This does not hold true for the third cohomology group $\mathrm{H}^3(\mathcal{V},\mathcal{V})$. In fact, in the next subsection we will give a sketch of the proof of the following results:
\begin{equation*}
dim(\mathrm{H}^3(\mathcal{V},\mathbb{K}))=1 \qquad\text{ and }\qquad dim(\mathrm{H}^3(\mathcal{V},\mathcal{V}))=1\,. 
\end{equation*} 
Hence, the third cohomology of the Virasoro algebra is not zero, but one-dimensional.
\subsection{The third cohomology of the Virasoro algebra}
\newtheorem{H3Virasoro2}{Theorem}[subsection]
\begin{H3Virasoro2} \label{h3K}
The third cohomology groups of the Witt algebra $\mathcal{W}$ and the Virasoro algebra $\mathcal{V}$ over a field $\mathbb{K}$ with $char(\mathbb{K})=0$ and values in the trivial module are one-dimensional, i.e.
\begin{equation*}
dim(\mathrm{H}^3(\mathcal{W},\mathbb{K}))=dim(\mathrm{H}^3(\mathcal{V},\mathbb{K}))=1\,.
\end{equation*}
\end{H3Virasoro2}
\begin{proof}
In a first step, we prove that the dimensions of $\mathrm{H}^3(\mathcal{W},\mathbb{K})$ and $\mathrm{H}^3(\mathcal{V},\mathbb{K})$ are bounded from above by one, i.e. $dim(\mathrm{H}^3(\mathcal{W},\mathbb{K}))\leq 1$ and $dim(\mathrm{H}^3(\mathcal{V},\mathbb{K}))\leq 1$. The proof of this result is almost the same for the Witt algebra and the Virasoro algebra such that both Lie algebras can be analyzed simultaneously. The techniques used in the proof are elementary algebraic manipulations, similar to those used to prove $\mathrm{H}^3(\mathcal{W},\mathcal{W})=\{0\}$, c.f. \cite{Ecker}.\\
In a second step, we prove that the dimensions of $\mathrm{H}^3(\mathcal{W},\mathbb{K})$ and $\mathrm{H}^3(\mathcal{V},\mathbb{K})$ are bounded from below by one, i.e. $dim(\mathrm{H}^3(\mathcal{W},\mathbb{K}))\geq 1$ and $dim(\mathrm{H}^3(\mathcal{V},\mathbb{K}))\geq 1$. This is done as follows.\par
The fact that the dimension of $\mathrm{H}^3(\mathcal{W},\mathbb{K})$ is one agrees with the result obtained from continuous cohomology for $Vect(S^1)$, i.e. $dim(\mathrm{H}^3(Vect(S^1),\mathbb{R}))=1$, see \cite{MR874337,MR0245035}. Inspired by this result, we consider the Godbillon-Vey $3$-cocycle $\mathscr{GV}$, which is a generator of the space $\mathrm{H}^3(Vect(S^1),\mathbb{R})$, see e.g. the book by Guieu and Roger \cite{MR2362888} and also the original literature \cite{MR0245035}. Expressing the cocycle $\mathscr{GV}$ in our basis and in our algebraic setting, we obtain:
\begin{align*}
&\mathscr{GV}(e_n,e_m,e_k)=\left\{
\begin{array}[h]{cl}
	0 & \text{ for }k\neq -(n+m)\\
	A\left|
	\begin{array}[h]{ccc}
	1 & 1 & 1\\
	n & m & -(n+m)\\
	n^2 & m^2 & (n+m)^2
	\end{array}\right|	& \text{ for }k=n+m
		 \end{array}\right. \\
		\Leftrightarrow &\ \mathscr{GV}(e_n,e_m,e_k)=\left\{
\begin{array}[h]{cl}
	0 & \text{ for }k\neq -(n+m)\\
	A(m-n)(2m+n)(m+2n)	& \text{ for }k=n+m
		 \end{array}\right. \,,
\end{align*}  
where $A$ is a non-vanishing constant. A direct verification shows that $\mathscr{GV}$ is a $3$-cocycle for $\mathrm{H}^3(\mathcal{W},\mathbb{K})$. Moreover, evaluating $\mathscr{GV}$ on the generators $e_1,e_0$ and $e_{-1}$ yields a value different from zero for $\mathscr{GV}$. However, every coboundary evaluated on $e_1,e_0$ and $e_{-1}$ gives zero. Thus, the $3$-cocycle $\mathscr{GV}$ is not a coboundary. 
By defining $\mathscr{GV}(e_n,e_m,t)=0$, the Godbillon-Vey cocycle can be trivially extended to a cochain of $\mathrm{H}^3(\mathcal{V},\mathbb{K})$. A similar direct verification shows that the extended $\mathscr{GV}$ is a non-trivial cocycle of $\mathrm{H}^3(\mathcal{V},\mathbb{K})$. This shows that $\mathrm{H}^3(\mathcal{W},\mathbb{K})$ and $\mathrm{H}^3(\mathcal{V},\mathbb{K})$ are at least one-dimensional.\\
Together with the first step we obtain the statement of the theorem. 
  \end{proof}
\newtheorem{remark}{Remark}[subsection]
\begin{remark} 
We showed that ${H}^3(\mathcal{W},\mathbb{K})$ and ${H}^3(\mathcal{V},\mathbb{K})$ are generated by the cocycle class of the algebraization of the Godbillon-Vey cocycle.
\end{remark}
\newtheorem{H3Virasoro3}[H3Virasoro2]{Theorem}
\begin{H3Virasoro3} 
The third cohomology group of the Virasoro algebra $\mathcal{V}$ over a field $\mathbb{K}$ with $char(\mathbb{K})=0$ and values in the adjoint module is one-dimensional, i.e.
\begin{equation*}
dim(\mathrm{H}^3(\mathcal{V},\mathcal{V}))=1\,.
\end{equation*}
\end{H3Virasoro3}
\begin{proof}
This proof is similar, although distinctly more intricate, to the one provided explicitly in the appendix to prove $\mathrm{H}^1(\mathcal{V},\mathcal{V})=\{0\}$. We will provide a sketch of the proof.\\
The following exact sequence of Lie algebras,
\begin{equation*}
0\longrightarrow \mathbb{K}\longrightarrow \mathcal{V}\longrightarrow\mathcal{W}\longrightarrow 0\,,
\end{equation*}
is also an exact sequence of $\mathcal{V}$-modules. Such sequences give rise to long exact sequences in cohomology. 
The relevant part for us is the part related to the third cohomology,
\begin{equation}
\dots \longrightarrow \mathrm{H}^2(\mathcal{V},\mathcal{W})\longrightarrow \mathrm{H}^3(\mathcal{V},\mathbb{K})\longrightarrow \mathrm{H}^3(\mathcal{V},\mathcal{V})\longrightarrow\mathrm{H}^3(\mathcal{V},\mathcal{W})\longrightarrow\dots\,.\label{LongSequence}
\end{equation}
In \cite{MR3200354}, it was shown that $\mathrm{H}^2(\mathcal{V},\mathcal{W})\cong\mathrm{H}^2(\mathcal{W},\mathcal{W})=\{0\}$. Similarly, by the same kind of computations as the ones presented in the main part of this article we show (see \cite{Ecker}) that:
\begin{equation*}
\mathrm{H}^3(\mathcal{V},\mathcal{W})\cong\mathrm{H}^3(\mathcal{W},\mathcal{W})\,,
\end{equation*}
and by taking into account Theorem \ref{h3}, we obtain:
\begin{equation}
\mathrm{H}^3(\mathcal{V},\mathcal{W})=\{0\}\,.\label{H3VW}
\end{equation}
The proof is an immediate generalization of the one performed in \cite{MR3200354}. \\
Taking into account the result (\ref{H3VW}) and $\mathrm{H}^2(\mathcal{W},\mathcal{W})=\{0\}$, the long exact sequence (\ref{LongSequence}) becomes:
\begin{equation*}
0\longrightarrow \mathrm{H}^3(\mathcal{V},\mathbb{K})\longrightarrow\mathrm{H}^3(\mathcal{V},\mathcal{V})\longrightarrow 0
\end{equation*}
By Theorem \ref{h3K}, we obtain $dim(\mathrm{H}^3(\mathcal{V},\mathcal{V}))=dim(\mathrm{H}^3(\mathcal{V},\mathbb{K}))=1$.
\end{proof}
\newtheorem{remark2}[remark]{Remark}\label{Remark2}
\begin{remark2} 
A shorter alternative proof of $\mathrm{H}^3(\mathcal{V},\mathcal{W})\cong\mathrm{H}^3(\mathcal{W},\mathcal{W})$ can be obtained via the Hochschild-Serre spectral sequence. In our case, the second page of the Hochschild-Serre spectral sequence corresponds to $E_2^{p,q}=\mathrm{H}^p(\mathcal{W},\mathrm{H}^q(\mathbb{K},M))$. Taking $\mathcal{W}$ as $\mathcal{V}$-module and denoting by $\varphi_i:\mathrm{H}^i(\mathcal{W},\mathcal{W})\rightarrow\mathrm{H}^{i+2}(\mathcal{W},\mathcal{W})$ the map on the $E_2$ level, we obtain from the third page $E_3=E_\infty$ the following result:
\begin{equation*}
\mathrm{H}^k(\mathcal{V},\mathcal{W})\cong\frac{\mathrm{H}^k(\mathcal{W},\mathcal{W})}{\text{im}\ \varphi_{k-2}}\bigoplus\text{ker}\left(\varphi_{k-1}:\mathrm{H}^{k-1}(\mathcal{W},\mathcal{W})\rightarrow \mathrm{H}^{k+1}(\mathcal{W},\mathcal{W})\right)\,,
\end{equation*}  
and in particular, if $\mathrm{H}^j(\mathcal{W},\mathcal{W})=\{0\}$  for  $k-1\leq j\leq k$, then 
\begin{equation*}
\mathrm{H}^{k+1}(\mathcal{V},\mathcal{W})\cong\mathrm{H}^{k+1}(\mathcal{W},\mathcal{W})\,.
\end{equation*}
Applying this to $k=2$ and using the results given in Section \ref{ResultsWitt}, we obtain $\mathrm{H}^3(\mathcal{V},\mathcal{W})\cong\mathrm{H}^3(\mathcal{W},\mathcal{W})$. 
\end{remark2}
\appendix
\section{The first cohomology}
In this appendix, we provide an explicit proof of the vanishing of the first cohomology groups, $\mathrm{H}^1(\mathcal{W},\mathcal{W})=\{0\}$ and $\mathrm{H}^1(\mathcal{V},\mathcal{V})=\{0\}$, as we want to present a gentle introduction to our techniques. 
\subsection{The Witt algebra}
In this section, we analyze the first cohomology group of the Witt algebra with values in the adjoint module. It is already known that all derivations of the Witt algebra are inner derivations, see Zhu and Meng \cite{MR1775845}, i.e. $\mathrm{H}^1(\mathcal{W},\mathcal{W})=0$ in the language of cohomological algebra. 
However, we shall prove this result again by using our algebraic techniques developed here, in order to introduce our notation and to provide a warm-up example of the used procedures.
The main aim of this section is to proof the following theorem:
\newtheorem{H1}{Theorem}[subsection]
\begin{H1} \label{h1}
The first cohomology of the Witt algebra $\mathcal{W}$ over a field $\mathbb{K}$ with char($\mathbb{K}$)$=0$ and values in the adjoint module vanishes, i.e.
\begin{equation*}
\mathrm{H}^1(\mathcal{W},\mathcal{W})=\{0\}\,.
\end{equation*}
\end{H1}
The proof follows in two steps, the first step concentrating on the non-zero degree part of the Witt algebra, the second step focusing on the degree zero part.\\
Recall that the coboundary condition for a 1-cocycle $\psi\in \mathrm{H}^1(\mathcal{W},\mathcal{W})$ is given by:
\begin{equation*}
\psi(x)=(\delta_0\phi)(x)=-x\cdot\phi=\left[\phi,x\right]\ ,
\end{equation*}
with $x\in\mathcal{W}$ and $\phi\in \mathrm{C}^0(\mathcal{W},\mathcal{W})=\mathcal{W}$. As the values are taken in the adjoint module, $\cdot=\left[.,.\right]$.\\
The cocycle condition for a 1-cocycle $\psi$ is given by:
\begin{equation*}
\delta_1\psi(x_1,x_2)=0=\psi\left(\left[x_1,x_2\right]\right)-\left[x_1,\psi(x_2)\right]-\left[\psi(x_1),x_2\right]\ ,
\end{equation*}
with $x_1,x_2\in\mathcal{W}$.
\subsubsection{The non-zero degree part of the Witt algebra}
 \newtheorem{H1DegreeNonZero}{Proposition}[H1]
\begin{H1DegreeNonZero} 
The following hold:
\begin{align*}
&\mathrm{H}^1_{(d)}(\mathcal{W},\mathcal{W})=\{0\}\ \text{ for }\ d\neq 0\ ,\\
&\mathrm{H}^1(\mathcal{W},\mathcal{W})=\mathrm{H}^1_{(0)}(\mathcal{W},\mathcal{W})\ .
\end{align*}
\end{H1DegreeNonZero}
\begin{proof} Let $\psi\in\mathrm{H}^1_{(d\neq 0)}(\mathcal{W},\mathcal{W})$. \\
Let us perform a cohomological change $\psi'=\psi -\delta_0\phi$ with the following 0-cochain $\phi$:
\begin{equation*}
\phi=-\frac{1}{d}\psi(e_0)\in\mathcal{W}\Rightarrow (\delta_0\phi)(x)=\frac{1}{d}\left[x,\psi(e_0)\right]\ ,
\end{equation*}
which gives us:
\begin{equation*}
\psi'(x)=\psi(x) -(\delta_0\phi)(x)=\psi(x)+\frac{1}{d}\left[\psi(e_0),x\right]\,.
\end{equation*}
Hence,
\begin{equation*}
\psi'(e_0)=\psi(e_0)+\frac{1}{d}\left[\psi(e_0),e_0\right]=\psi(e_0)-\frac{1}{d}\ \text{deg}(\psi(e_0))\ \psi(e_0)
=\psi(e_0)-\frac{1}{d}\ d\ \psi(e_0)=0\,.
\end{equation*}
We thus have \framebox{$\psi'(e_0)=0$}.\\
Next, let us write down the cocycle condition for $\psi'$ on the doublet $(x,e_0)$ for $x\in\mathcal{W}$:
\begin{align*}
& 0=\psi'\left(\left[x,e_0\right]\right)-\underbrace{\left[x,\psi'(e_0)\right]}_{=0}-\left[\psi'(x),e_0\right]
\Leftrightarrow\   0=\psi'(-\text{deg}(x)x)+\text{deg}(\psi'(x))\psi'(x)\\
&\Leftrightarrow\   0=-\text{deg}(x)\psi'(x)+(\text{deg}(x)+d)\psi'(x)
\Leftrightarrow\   0=d\ \psi'(x)\,.
\end{align*}
As $d\neq 0$, we get $\psi'(x)=0\ \forall x\in\mathcal{W}$, meaning that $\psi$ is a coboundary on $\mathcal{W}$.
We conclude that the first cohomology of the Witt algebra reduces to the degree zero part, in accordance with the result (\ref{Fuks}).
\end{proof}
\subsubsection{The degree zero part for the Witt algebra}
 \newtheorem{H1DegreeZero}{Proposition}[subsubsection]
\begin{H1DegreeZero} 
The following holds:
\begin{equation*}
\mathrm{H}^1_{(0)}(\mathcal{W},\mathcal{W})=\{0\}\ .
\end{equation*}
\end{H1DegreeZero}
\begin{proof}
Let $\psi$ be a degree zero 1-cocycle, i.e. we can write it as $\psi(e_i)=\psi_ie_i$ with suitable coefficients $\psi_i\in\mathbb{K}$. Consider the following 0-cochain $\phi=\psi_1 e_0$. The coboundary condition for $\phi$ gives:
\begin{equation*}
(\delta_0\phi)(e_i)= \left[\phi,e_i \right]=i\psi_1e_i\,.
\end{equation*}
The cohomological change $\psi'=\psi-\delta_0\phi$ leads to $\psi'_1=0$. In the following, we will work with a 1-cocycle normalized to $\psi_1'=0$, although we will drop the apostrophe in order to augment readability. \\
The 1-cocycle condition for $\psi$ on the doublet $(e_i,e_j)$ becomes:
\begin{align*}
& 0=\psi\left(\left[e_i,e_j\right]\right)-\left[e_i,\psi(e_j)\right]-\left[\psi(e_i),e_j\right]\\
\Leftrightarrow\  & 0=(j-i)\left(\psi_{i+j}-\psi_{j}-\psi_{i}\right)\,.
\end{align*}
For $j=1$ and $i=0$, we obtain from the 1-cocycle condition: $\psi_0=0$.\\
For $j=1$ and $i<0$ decreasing, we obtain from the 1-cocycle condition: $\psi_{i}=\psi_{i+1}=0$.\\
For $j=1$ and $i>1$ increasing, we obtain from the 1-cocycle condition: $\psi_{i+1}=\psi_{i}=\psi_{2}$, where the value of $\psi_2$ is unknown for the moment.\\
Next, taking $j=2$ and for example $i=3$, we obtain: 
\begin{align*}
&\psi_5-\psi_2-\psi_3=0\\
\Leftrightarrow\  &\psi_2-\psi_2-\psi_2=0\text{ as we have }\psi_{i}=\psi_{2}\ \forall i>1\\
\Leftrightarrow\  &\psi_2=0\,.
\end{align*}
All in all, we conclude \framebox{$\psi_i=0\ \forall i\in\mathbb{Z}$}
\end{proof}
This concludes the proof of Theorem \ref{h1}.

\subsection{The Virasoro algebra}
Recall from Sections \ref{ResultsWitt} and \ref{ResVirasoro} that we have $\mathrm{H}^1(\mathcal{W},\mathcal{W})~=~\{0\}$ and $\mathrm{H}^1(\mathcal{V},\mathbb{K})~=~\{0\}$, respectively. These results will be used in the proof of $\mathrm{H}^1(\mathcal{V},\mathcal{V})~=~\{0\}$ based on long exact sequences. This time we will focus only on the degree zero part accordingly to the result (\ref{Fuks}).
\newtheorem{H1Virasoro2}{Theorem}[subsection]
\begin{H1Virasoro2} 
The first cohomology group of the Virasoro algebra with values in the adjoint module is zero, i.e.:
\begin{equation*}
\mathrm{H}^1(\mathcal{V},\mathcal{V})=\{0\}\,.
\end{equation*}
\end{H1Virasoro2}
\begin{proof}
For the first cohomology, the relevant part of the long exact sequence (\ref{LongSequence}) looks as follows:
\begin{equation*}
\dots \longrightarrow \mathrm{H}^1(\mathcal{V},\mathbb{K})\longrightarrow \mathrm{H}^1(\mathcal{V},\mathcal{V})\stackrel{\nu_*}{\longrightarrow}\mathrm{H}^1(\mathcal{V},\mathcal{W})\longrightarrow\dots \,,
\end{equation*}
where $\nu:\mathcal{V}\rightarrow \mathcal{W}$ is the projection.
Since we already have $\mathrm{H}^1(\mathcal{V},\mathbb{K})=\{0\}$ and also $\mathrm{H}^1(\mathcal{W},\mathcal{W})=\{0\}$, it suffices to prove $\mathrm{H}^1(\mathcal{V},\mathcal{W})\cong \mathrm{H}^1(\mathcal{W},\mathcal{W})$ in order to conclude. Using the result from Remark \ref{Remark2} for $k=0$, and recalling from Section \ref{ResultsWitt} that $\mathrm{H}^0(\mathcal{W},\mathcal{W})=0$, this follows immediately. Nevertheless, as Remark \ref{Remark2} is based on the Hochschild-Serre spectral sequence, which is contrary to our intention to stay elementary, we prefer to present a more direct proof in the following.\\ 
The proof consists of two steps. First, we will compare the cocycles of $Z^1(\mathcal{V},\mathcal{W})$ to the cocycles of $Z^1(\mathcal{W},\mathcal{W})$. In the second step, we will compare the coboundaries of $B^1(\mathcal{V},\mathcal{W})$ to the ones of $B^1(\mathcal{W},\mathcal{W})$.\\
Let $\hat{\psi}:\mathcal{V}\rightarrow\mathcal{W}$ be a cocycle of $Z^1(\mathcal{V},\mathcal{W})$. Our aim is to show that the restriction of this cocycle to $\mathcal{W}$, i.e. $\psi:=\hat{\psi}|_\mathcal{W}:\mathcal{W}\rightarrow\mathcal{W}$, is a cocycle of $Z^1(\mathcal{W},\mathcal{W})$. Let $x_1,x_2\in\mathcal{V}$. Writing the Virasoro product $[\cdot,\cdot]_\mathcal{V}$ in terms of the Witt product $[\cdot,\cdot]_\mathcal{W}$ and the 2-cocycle $\alpha(\cdot,\cdot)$ giving the central extension, i.e. $[\cdot,\cdot]_\mathcal{V}=[\cdot,\cdot]_\mathcal{W}+\alpha(\cdot,\cdot)\cdot t$, the cocycle condition for $\hat{\psi}$ becomes:
\begin{align}
&0=(\delta_1^\mathcal{V}\hat{\psi})(x_1,x_2)=\hat{\psi}\left([x_1,x_2]^\mathcal{V}\right)-x_1\cdot\hat{\psi}(x_2)+x_2\cdot\hat{\psi}(x_1)\nonumber\\
&\Leftrightarrow 0=(\delta_1^\mathcal{V}\hat{\psi})(x_1,x_2)=\hat{\psi}\left([x_1,x_2]^\mathcal{W}\right)+\alpha(x_1,x_2)\ \hat{\psi}(t)-[x_1,\hat{\psi}(x_2)]^\mathcal{W}+[x_2,\hat{\psi}(x_1)]^\mathcal{W}\nonumber\\
&\Leftrightarrow 0=(\delta_1^\mathcal{V}\hat{\psi})(x_1,x_2)=(\delta_1^\mathcal{W}\hat{\psi})(x_1,x_2)+\alpha(x_1,x_2)\ \hat{\psi}(t)\,.\label{CompCoc}
\end{align} 
 Since we are considering degree-zero cocycles, the cocycle $\hat{\psi}$ evaluated on the central element reads as follows:
\begin{equation*}
\hat{\psi}(t)=c\ e_0\,,
\end{equation*}
for suitable $c\in\mathbb{K}$.
Next, let us insert this expression into the cocycle condition for $(e_1,t)$, which yields:
\begin{align*}
&(\delta_1^\mathcal{V}\hat{\psi})(e_1,t)=\hat{\psi}\left(\cancel{[e_1,t]^\mathcal{V}}\right)-e_1\cdot\hat{\psi}(t)+\cancel{t\cdot\hat{\psi}(e_1)}=0\\
&\Leftrightarrow -[e_1,\hat{\psi}(t)]^\mathcal{W}=-c\ [e_1,e_0]^\mathcal{W}=c\ e_1=0\\
&\Leftrightarrow c=0\,.
\end{align*}
 Inserting $\hat{\psi}(t)=0$ into (\ref{CompCoc}), we obtain
\begin{equation}
0=(\delta_1^\mathcal{V}\hat{\psi})(x_1,x_2)=(\delta_1^\mathcal{W}\psi)(x_1,x_2)\,.\label{CompCoc2}
\end{equation}
This means that a cocycle $\hat{\psi}\in Z^1(\mathcal{V},\mathcal{W})$ corresponds to a cocycle $\psi\in Z^1(\mathcal{W},\mathcal{W})$ when projected to $\mathcal{W}$. Moreover, a cocycle $\psi$ of $Z^1(\mathcal{W},\mathcal{W})$ can also be lifted to a cocycle $\hat{\psi}:=\psi\circ\nu$ in $Z^1(\mathcal{V},\mathcal{W})$. By definition, we thus have $\hat{\psi}(t)=0$ and the relation (\ref{CompCoc2}) holds true. Hence, a cocycle $\psi\in Z^1(\mathcal{W},\mathcal{W})$ yields a cocycle $\hat{\psi}\in Z^1(\mathcal{V},\mathcal{W})$ and we have $Z^1(\mathcal{V},\mathcal{W})\cong Z^1(\mathcal{W},\mathcal{W})$ in a canonical way.\\
The second step of the proof consists in comparing the coboundaries of $B^1(\mathcal{V},\mathcal{W})$ and those of $B^1(\mathcal{W},\mathcal{W})$. However, this is trivial. In fact, the coboundary condition applied on a 0-cochain $\phi\in C^0(\mathcal{V},\mathcal{W})$ is the same as the one applied on a 0-cochain $\phi\in C^0(\mathcal{W},\mathcal{W})$, yielding in both cases:
\begin{equation*}
(\delta_0\phi)(x)=-x\cdot\phi\text{ with }\phi\in\mathcal{W}\,.
\end{equation*}
Since the central element of $\mathcal{V}$ acts trivially on $\mathcal{W}$, we have $B^1(\mathcal{V},\mathcal{W})\cong B^1(\mathcal{W},\mathcal{W})$. All in all, we conclude $\mathrm{H}^1(\mathcal{V},\mathcal{W})\cong \mathrm{H}^1(\mathcal{W},\mathcal{W})$ in a canonical way.
\end{proof}
\bibliographystyle{amsplain}
\bibliography{refs_cohom_lie}

\end{document}